\documentclass[reqno]{article}
\usepackage{authblk}
\pdfoutput=1
\usepackage{amsmath}
\usepackage{amsfonts,amsthm,amssymb}
\usepackage{mathabx}
\usepackage{bm}
\usepackage{amscd}
\usepackage{comment}
\usepackage[latin2]{inputenc}
\usepackage{t1enc}
\usepackage[mathscr]{eucal}
\usepackage{indentfirst}
\usepackage{graphicx}
\usepackage{graphics}
\usepackage{pict2e}
\usepackage{epic}
\numberwithin{equation}{section}
\usepackage[margin=3cm]{geometry}
\usepackage{epstopdf} 
\usepackage[colorlinks,linkcolor=blue]{hyperref}
\usepackage[capitalise,noabbrev]{cleveref}

\usepackage{todonotes}
\crefformat{equation}{(#2#1#3)}
\crefmultiformat{equation}{(#2#1#3)}{ and~(#2#1#3)}{, (#2#1#3)}{ and~(#2#1#3)}
\crefrangeformat{equation}{(#3#1#4) to~(#5#2#6)}
\usepackage[normalem]{ulem}

\allowdisplaybreaks

\theoremstyle{plain}
\newtheorem{Thm}{Theorem}[section]
\newtheorem*{Thm*}{Theorem}
\newtheorem{Lem}[Thm]{Lem}

\theoremstyle{definition}

\newtheorem{Rem}[Thm]{Rem}
\newtheorem{?}[Thm]{Problem}

\newcommand{\R}{\mathbb{R}}

\begin{document}
	
	\begin{titlepage}
		\title{Asymptotic stability of the composite wave of rarefaction wave and contact wave to nonlinear viscoelasticity model with non-convex flux}
        \author{Zhenhua Guo$^{1}$} 
	\author{Meichen Hou$^{2}$\thanks{ Corresponding author.}}
	\author{Guiqin Qiu$^{3}$}
		\author{Lingda Xu$^{4}$}



		\affil{
            \footnotesize\qquad\quad $^1$ School of Mathematics and Information Science, Guangxi University, Nanning 530004, P.R.China.\\
           E-mails: zhguo@gxu.edu.cn}
            \affil{\footnotesize\qquad\quad $^2$  Center for nonlinear studies, School of Mathematics, Northwest University, Xi'an 710069, P.R.China.\\
            E-mails: meichenhou@nwu.edu.cn}
	     
	      \affil{\footnotesize\qquad\quad $^3$  Academy of Mathematics and Systems Science, Chinese Academy of Sciences, Beijing 100190, P.R.China.\\
	   E-mails: guiqinqiu@126.com
    }
            \affil{\footnotesize\qquad\quad $^4 $ Department of Applied Mathematics,
The Hong Kong Polytechnic University,
Hong Kong,
P.R.China.\\
E-mails: lingda.xu@polyu.edu.hk}
		\date{}

	\end{titlepage}
	
	\maketitle

\begin{abstract}

In this paper, we consider the wave propagations of viscoelastic materials, which has been derived by Taiping-Liu to approximate the viscoelastic dynamic system with fading memory (see [T.P.Liu(1988)\cite{LiuTP}]) by the Chapman-Enskog expansion. By constructing a set of linear diffusion waves coupled with the high-order diffusion waves to achieve cancellations to approximate the viscous contact wave well and explicit expressions, the nonlinear stability of the composite wave is obtained by a continuum argument.

It emphasis that,
the stress function in our paper is a general non-convex function, which leads to several essential differences from strictly hyperbolic systems such as the Euler system.  Our method is completely new and can be applied to more general systems and a new weighted Poincar\'e type of inequality is established, which is more challenging compared to the convex case and this inequality plays an important role in studying systems with non-convex flux.

	\bigbreak

	
\end{abstract}

\tableofcontents

\section{Introduction }

\subsection{The nonlinear system of viscoelasticity}

In this paper, we study the large-time behavior of solutions to the nonlinear system of viscoelasticity. Such problems arise naturally in conducting physical experiments on viscoelastic materials and manufacturing of polymeric chemical materials.

Our focus is on the viscoelastic dynamic system with fading memory mentioned in T.-P. Liu's paper \cite{LiuTP}. Roughly speaking,  the motion of an unbounded, homogeneous, viscoelastic bar with fading memory is governed by the following equations:
\begin{align}\label{eq1}
\left\{\begin{aligned}
&v_t- u_x  =0,\\
&u_t- \tilde{\sigma}(v)_x =0,\quad x \in \mathbb{R}, \quad t>0,
\end{aligned}\right.
\end{align}
where $v, u$, and $\tilde{\sigma}$ are the strain deformation, velocity, and stress; the stress $\tilde{\sigma}(x,t)$ is a given function of the strain $v$ and its past history,
\begin{align}\label{sigmag}
\begin{aligned}
\tilde{\sigma}(x, t)=f(v(x, t))+\int_{-\infty}^t a^{\prime}(t-\tau) g(v(x, \tau)) d \tau ;
\end{aligned}
\end{align}
here $a(s)$ is a given kernel on $0 \leqq s<\infty$ with derivative $a^{\prime}$, and $f(v), g(v)$ are given smooth material functions. The history and initial data are given by
$$
v(x, t)=\eta(x, t), \quad t \leqq 0, \quad u(x, 0)=u_0(x) .
$$

As in \cite{LiuTP}, by the Chapman-Enskog expansion,  the hyperbolic-parabolic system
\begin{align}\label{sys1}
\left\{\begin{aligned}
&v_t-u_x=0 \\
&u_t-\sigma(v)_x=\left(\mu(v) u_x\right)_x,
\end{aligned}\right.
\end{align}
has been proposed to approximate \eqref{eq1}, where

$$
\begin{aligned}
& \sigma(v)=f(v)-a(0) g(v), \qquad\qquad \mu(v)=g^{\prime}(v) \int_{-\infty}^0 s a^{\prime}(-s) d s>0.
\end{aligned}
$$
Here we denote the approximate stress function of strain as $\sigma(v)$ when compared with $\tilde{\sigma}$ in \eqref{eq1}.

There exists large literatures on experiments and numerical simulations about the viscoelastic materials. However, a rigorous mathematical analysis theory is still being completed \cite{CoGuHe,Da-No,GR,Pi,RHN}.


 Initially, scholars established corresponding mathematical models for different types of viscoelastic materials, such as Kelvin-Voigt material, Maxwell material,  standard linear solid model, etc. see \cite{MBB,RD,TR} and the references therein. In recent years, the emergence of polymer materials and composite materials, as well as the need to analyze the mechanical response and stability of some major engineering structures, has prompted people to pay sufficient attention to viscoelastic theory. From the end of the 1950s to the early 1960s, Coleman and Noll et al developed the constitutive theory of viscoelastic materials with memory. Especially, the article "Thermodynamic of materials with memory" \cite{Cole}, given by Coleman in 1964, exerted a great influence on this domain. Since then, many scholars have conducted in-depth research on viscoelastic dynamic systems with memory.
In 1976, Greenberg-Hastings\cite{GR-HA} used the monotonic methods to prove the existence of traveling wave solutions to the equations of motions for materials that exhibit long-range memory effects. Subsequently, in 1988, Liu\cite{LiuTP} proved the existence of smooth and non-smooth traveling wave solutions for one-dimensional nonlinear viscoelastic dynamic equations in materials with fading memory. As for the high-dimensional case, Qin-Ni \cite{Qin-Ni} proposed a high-order iterative method and they proved the existence of traveling wave solutions for the three-dimensional nonlinear viscoelastic dynamic equations with special integral kernels. All the aforementioned studies have discovered the phenomenon of wave propagation in viscoelastic media. In addition, there are also some discussions on the existence of global solutions and periodic solutions for nonlinear viscoelastic dynamic equations. We refer to \cite{BN,FE,Hr-No,Macc,Staffans}.


\subsection{The wave propagation in the viscoelastic system and its differences compared with the Euler system.}

We are interested in studying the wave propagation for \eqref{sys1} under the assumption that
the stress function $\sigma(v)$ is non-convex, i.e. $\sigma''(v)>0$ doesn't hold in $\R$ and $\sigma''(v)\equiv 0$ on some interval $[c,d]$, see \eqref{sigma-1} below. Now we want to consider the Cauchy problem of \eqref{sys1}. For simplicity, we assume that  $\mu>0$ be a constant as \cite{KM1994,Liu-shock,MM}, that is,
\begin{align}\label{eq}
\left\{\begin{aligned}
&v_t- u_x  =0,\\
&u_t- \sigma(v)_x =\mu u_{xx}, \ \ \ (x,t)\in \mathbb{R}\times \mathbb{R}_{+}.
\end{aligned}\right.
\end{align}
 Initial condition of \eqref{eq} satisfies
\begin{equation}\label{inta}
\begin{aligned}
(v,u)(x,0)=(v_0,u_0)(x) \rightarrow (v_{\pm},u_{ \pm}), \quad \text { as } x \rightarrow \pm \infty.
\end{aligned}
\end{equation}
The far field states $(v_{\pm},u_{\pm})$ are constants,  we assume that $v_{-}<v_{+}$ without loss of generality.
It is well known that asymptotic behavior for \eqref{eq} is governed by the Riemann solutions of following inviscid system, i.e.,
\begin{align}\label{non-viscosity}
\left\{\begin{array}{l}
v_t-u_x  =0,\\[2mm]
u_t-\sigma(v)_x =0, \ \ \ (x,t)\in \mathbb{R}\times \mathbb{R}_{+}.
\end{array}\right.
\end{align}
If $v>0$, the system \eqref{non-viscosity} is the famous Euler system, which is a typical strictly hyperbolic system. The main feature of the Euler system is the singularity formation, that is, no matter how smooth or small the initial data is, the shock may form in a finite time. Thus, the study of singular solutions is important for this kind of hyperbolic system. It is well-known that the two characteristic fields of the Euler system are both generally non-linear, thus shock wave and rarefaction wave will form.

The emergence of these wave phenomena is an important feature of hyperbolic systems of conservation laws. Thus, it is of great importance to study the large-time behavior of solutions of the corresponding viscous systems, such as \eqref{eq}.  Hence, large-time behavior of solutions for system \eqref{eq} whose initial data and histories are smooth is studied by many scholars.  There is fruitful research about the stability of shock and rarefaction waves for \eqref{eq} with $v>0$, we refer to \cite{KM1985,MN1,MN2}.

The problem we consider in this paper is more complicated and has essential differences from the Euler system. Although the form of the equation is the same, the crucial condition of $v>0$ no longer holds in the viscoelastic model. That means \eqref{non-viscosity} is not a strictly hyperbolic system, there are much more rich wave phenomena. The Riemann solutions of the viscoelastic system \eqref{non-viscosity} consist of elementary waves such as shock waves, rarefaction waves, and contact discontinuities and their superposition,
which is essentially different from the Euler system. Moreover, only one wave pattern will be generated in one characteristic field in the case of the Euler system and different wave patterns in different characteristic fields will separate from each other as time evolves. But in the case of the viscoelastic system, two or more kinds of different wave patterns will occur in the same characteristic field. Furthermore, these wave patterns generated in the same characteristics field will not separate from each other even the time goes to infinity, which is very interesting but difficult to study its behavior.

When the stress $\sigma(v)$ is non-convex, stability analysis of solutions for \eqref{eq}  is firstly studied in\cite{KM1994,MM}. The stress function $\sigma(v)\in C^{1}(\mathbb{R})$ in\cite{KM1994,MM} satisfies  one of the following conditions  respectively :
\begin{equation}\label{sigma}
\left\{
\begin{aligned}
&\sigma'(v)>0,\quad  \sigma''(v)>0, v>0,\quad \sigma''(v)<0,v<0.\\
&\sigma'(v)>0,\quad \sigma''(v)<0, v>0,\quad \sigma''(v)>0, v<0.
\end{aligned}
\right.
\end{equation}
Precisely to say, Kawashima-Matsumura \cite{KM1994} proved the asymptotic stability of shock profiles for system \eqref{eq} under assumption $\eqref{sigma}_1$. Later on,
Matsumura-Mei \cite{MM} consider another case $\eqref{sigma}_2,$  and they also obtain the stability results of a shock wave for system \eqref{eq}. Both in \cite{KM1994,MM}, some delicate weight functions are introduced in the energy estimates to overcome the difficulty caused by the non-convexity of the characteristic field.
It can be observed from \eqref{sigma} that $v=0$ is an inflection point of concavity and convexity changes of the function $\sigma(v)$. 


Also, Nishihara proved the stability of degenerate shock for $\eqref{eq}$-$\eqref{sigma}_1$. Here are some other extended results, see \cite{AG,MN1994,MeiNishi}. Recently, Liu-Guo\cite{LG} gave the stability analysis of the large amplitude viscous shock waves for system \eqref{eq}. In \cite{LG}, They introduce an effective velocity to deform system \eqref{eq} so that the weight functions are no longer required in energy estimates.  Note that when $\sigma''(v)<0,v\in\mathbb{R}$, system \eqref{eq} is reduced to the $p$-system, which describes the motion of compressible viscous gas. The research on $p$-system has already been mature, we refer to \cite{KM1985,MN1,JS, KVW1} and the references therein. For general results about the global existence of solutions and their stability analysis for \eqref{eq1}, we refer to \cite{Shu-Zeng,Zeng}.

Since there have been remarkable results in the case of single-wave, the goal of this paper is to study the multi-wave patterns generated in the same characteristic field. Motivated by the research of \cite{MY} in single conservation law, which showed contact waves and rarefaction waves will generated if the flux is non-strictly-convex, we want to study the wave phenomenon for \eqref{eq} with general non-convex flux.

\subsection{Applications to physical model: The typical one-dimensional elastoplastic loaded wave propagation.}\label{physical}

In the mechanics of materials, the basic characteristics of one-dimensional elastoplastic loaded wave propagation in several typical elastoplastic materials are as follows:
\begin{itemize}
    \item Linear elastic materials only propagate elastic waves, and the wave speed and waveform remain unchanged, see Figure \ref{fig1} for example.
     \begin{figure}[h]
    \centering
\includegraphics[width=0.4\linewidth]{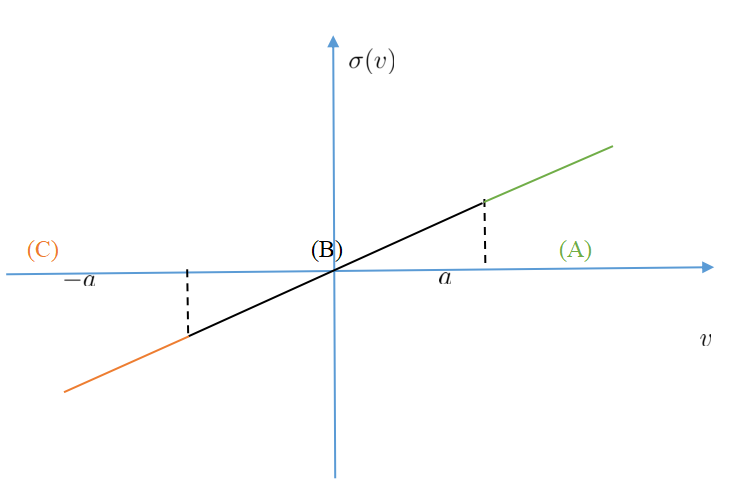}
    \caption{Linear elastic materials}
    \label{fig1}
\end{figure}
    \item Linear elasticity-linear hardening materials: two kinds of wave patterns propagate. One is an elastic wave with constant wave speed and waveform, but faster speed, as well as 
    plastic wave with  constant wave speed and waveform,  but slower speed.  The distance between the two waves continues to increase over time, see Figure \ref{fig2} for example.

     \begin{figure}[h]
    \centering
    \includegraphics[width=0.38\linewidth]{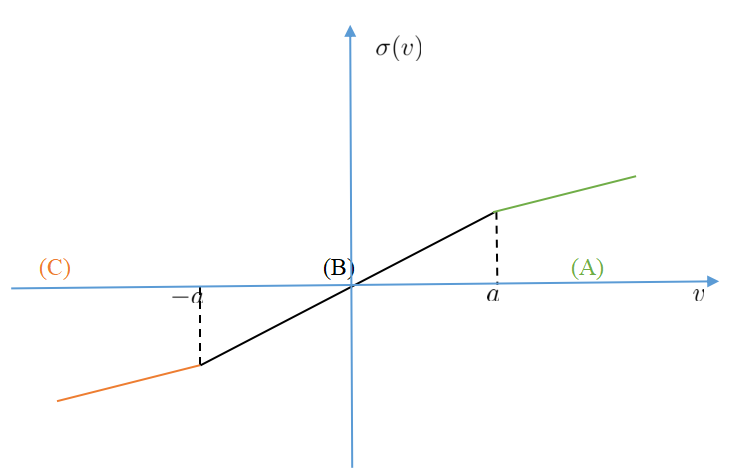}
    \caption{Linear elasticity-linear hardening materials}
    \label{fig2}
\end{figure}

    \item Linear elasticity-increasing hardening materials:  two kinds of elastoplastic waves propagate: include elastic waves with constant wave speed and waveform, but faster speed, as well as plastic waves whose wave speed and waveform change.  Plastic waves with larger amplitude have faster speed.  As time goes by, plastic waves with larger amplitude behind them will catch up with plastic waves with smaller amplitude before them, and the wave profile will become steeper and steeper.  Finally, convergence will occur, and the particle, velocity, and stress strain on the wavefront will produce a sudden jump, forming a shock wave.  The wavefront of the shock wave is a singular surface that must satisfy the shock jump conditions. See Figure \ref{fig3} for example.

     \begin{figure}[h]
    \centering
    \includegraphics[width=0.38\linewidth]{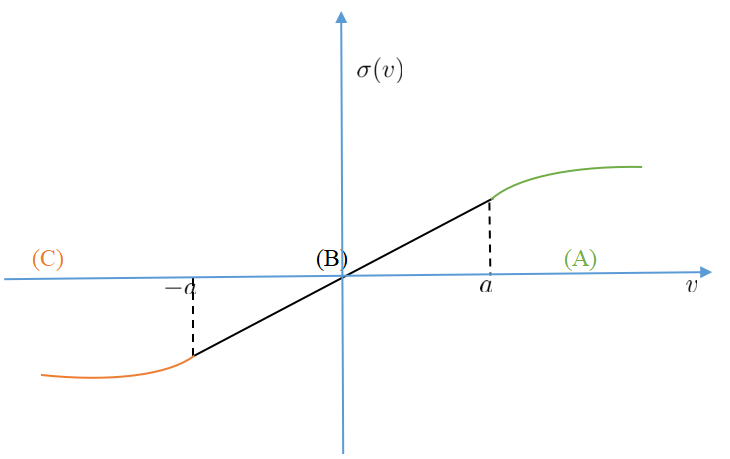}
    \caption{Linear elasticity-increasing hardening materials}
    \label{fig3}
\end{figure}

     \item Linear elasticity-decreasing hardening materials: there are two kinds of elastoplastic waves that propagate: elastic waves with constant wave speed and waveform, but faster speed, as well as plastic waves with changing wave speed and waveform, plastic waves with larger amplitude have slower speed, and the waveform of plastic waves will gradually widen and diverge as time goes by. See Figure \ref{fig4} for example.
     \begin{figure}[h]
    \centering
    \includegraphics[width=0.38\linewidth]{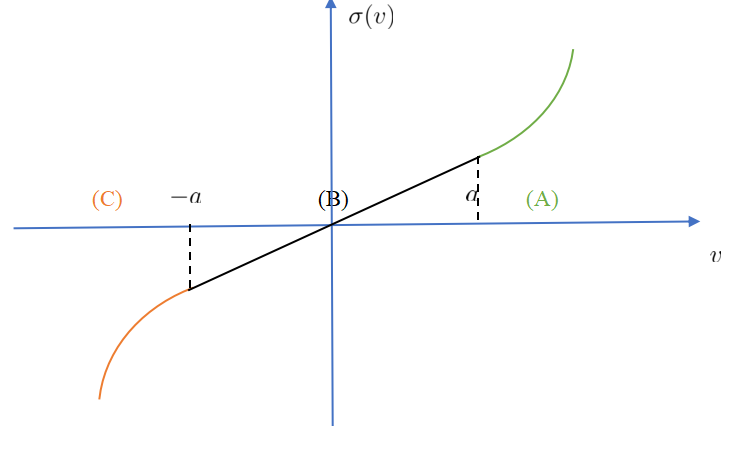}
    \caption{Linear elasticity-decreasing hardening materials}
    \label{fig4}
\end{figure}
\end{itemize}

When two elastoplastic loading waves propagating towards each other encounter each other, they will interact. When two loading waves with the same sign (both tensile or compression waves) encountering each other, they will be further loaded.  The complexity of the problem is also that for nonlinear problems, the superposition principle is no longer applicable.  The elastoplastic loaded wave will be reflected at the fixed end, and the elastoplastic wave will be gradually enhanced while reflecting back and forth between the fixed end and the impact end. And in this paper, we pay attention to study the wave propagation for the model of linear elasticity-decreasing hardening materials(see figure \ref{fig4}) , our stability results explain the wave interaction mechanics between elastic waves and plastic waves.

\subsection{Two waves propagate in the same characteristic field: difficulties and innovations.}
As explained before, we are interested in the phenomenon that two or more different wave patterns propagate in the same characteristic field. That is, we study the fourth case in subsection \ref{physical}. We describe the conditions of $\sigma(v)$ as follows:
\begin{equation}\label{sigma-1}
\left\{
\begin{aligned}
&\sigma''(v)>0,v\in[a,+\infty)=A,\\
&\sigma(v)=bv,b>0,v\in(-a,a)=B,\\
&\sigma''(v)<0,v\in(-\infty,-a]=C.
\end{aligned}
\right.
\end{equation}
where $b>0,a>0$ are positive constants, and  $\sigma'(v)$ is continuous on $(-\infty,+\infty),$ see figure \ref{fig4}.

Note that under assumption \eqref{sigma-1}, the characteristic field has both nonlinear and linearly degenerate components. We want to know how the propagation mechanism of waves will change in this general non-convex flow field. As for the Riemann problem of \eqref{non-viscosity}-\eqref{sigma-1}, different elementary waves on the same characteristic field may appear simultaneously. In particular, shock waves, rarefaction waves, and contact discontinuities may interact with each other within the same characteristic field. In this paper, we consider the case of a composite wave of the contact wave and rarefaction wave.
As the following picture, these two wave patterns will propagate in the same characteristic field.

    \begin{figure}[h]
    \centering
    \includegraphics[width=0.5\linewidth]{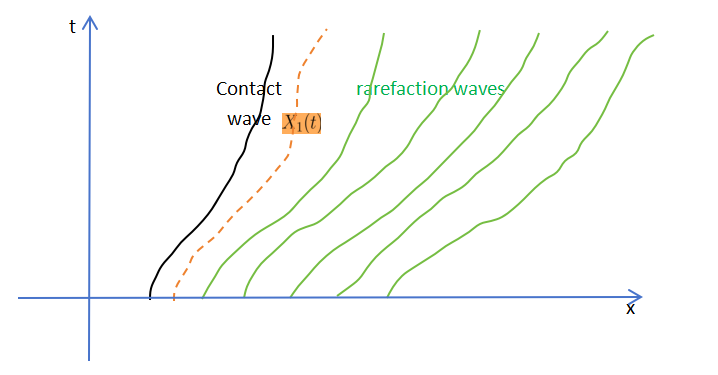}
    \caption{viscous contact wave and rarefaction wave propagate in the 2-characteristic field}
    \label{fig5}
   \end{figure}

Since the construction of wave patterns is complex we should explain it clearly, we give the rough version of the main result.
\begin{flushleft}
    \textbf{Rough version of the main result.}  The solution of viscoelastic system \cref{eq,inta} with the assumption \eqref{sigma-1} converges to the combination wave consisting of contact wave and rarefaction wave as time goes to infinity.
\end{flushleft}

This kind of nature makes the mathematical analysis difficult. Precisely,
\begin{itemize}
    \item The encounter of different waves may cause them to overlap and adhere to each other, resulting in a highly complex interaction mechanism that poses difficulty for our analysis.

    To overcome this difficulty, we introduce the adherent surface $X_1(t)$ to analyze the interaction of waves in the case of two overlapping viscous waves, see lemma \ref{xl}. This idea is inspired by Matsumura-Yoshida's stability analysis of nonlinear waves for non-convex scalar equation \cite{MY}.

    \item The construction of viscous contact wave.

    Another innovation lies in the construction of viscous contact waves.  We diagonalize system \eqref{eq} and introduce new variables $\Xi_i(i=1,2)$ in \eqref{Xi} to give an "approximate" solution for \eqref{eq}. Note that every component $\Xi_i$ satisfies a linear heat conduction equation and its explicit solution could be given by the heat-kernel function. Through this method, a certain combination of linear dissipative waves can be used to appropriately approximate the viscous contact discontinuities in our system.
    Different from the neat construction of the viscous contact wave in Huang-Matsumura-Xin \cite{HMX}, in which the contact wave of the Navier-Stokes equation is constructed based on a nonlinear diffusion wave, our contact wave is a combination of linear diffusion waves. More explicit expression of our viscous contact wave could be obtained and this approach can be applied to general viscous conservation law systems.

    \item A new Poincar\'e type of inequality.

    Note that due to the appearance of the viscous contact wave, there will be critical terms in the stability analysis
    \begin{align}
        \int_0^t\int_{\R}(1+\tau)^{-1}e^{-\frac{c\eta^2}{1+\tau}}(\phi^2,\psi^2)d\eta d\tau.
    \end{align}
    Inspired by Huang-Li-Matsumura \cite{HLM}, which introduced a key estimate on heat kernel, we developed a Poincar\'e type of inequality in lemma \ref{4.2} to estimate the critical terms.

    \item The delicate energy estimate.

    Due to the non-convexity of the flow field, we need to analyze the state of our solution in different regions of the flow field. A delicate estimate is constructed in lemma \ref{4.2}. At this point, the study on asymptotic behavior of solutions for corresponding viscous system \eqref{eq} is more challenging. We want to consider the large time behavior of \eqref{eq} under assumption \eqref{sigma-1}.
\end{itemize}

 This paper is arranged as follows: in section 2, we pay attention to analyzing the different wave patterns generated on the 2-characteristic fields for \eqref{eq},\eqref{sigma-1}, see subsection \ref{Sec2.1}. Then the construction of the viscous contact wave and its superposition with the rarefaction wave will be completed, see subsection \ref{Sec2.2}. We left the main theorem \ref{mt} at the end of this section. In section 3, we do the stability analysis to prove theorem \ref{mt}.

\section{The construction of composite wave patterns and the main results}

In this section, we turn to analyze the multi-wave patterns for system \eqref{eq} under assumptions \eqref{sigma-1}. Then we construct the ansatz and state the stability results in theorem \ref{mt}. 

\subsection{The Riemann solution for non-viscous system \eqref{non-viscosity}}\label{Sec2.1}
It is well known that the corresponding non-viscous system \eqref{non-viscosity} has two characteristic fields. The Jacobian matrix of flux function for \eqref{non-viscosity} is
\begin{align}
    \left(\begin{array}{cc}
        0&-1\\
        -\sigma(v)&0
    \end{array}\right),
\end{align}
and the corresponding eigenvalues are 

\begin{equation}\label{ev}
\begin{aligned}
\lambda_1=-\sqrt{\sigma'(v)},\quad \lambda_2=\sqrt{\sigma'(v)}.
\end{aligned}
\end{equation}
Moreover, 
\begin{equation}\label{sg}
\begin{aligned}
\nabla\lambda_i(v)\cdot r_i(v)=\frac{\sigma''(v)}{2\sqrt{\sigma'(v)}},
\end{aligned}
\end{equation}
where
\begin{equation}
\begin{aligned}
r_1(v)=(1,\sqrt{\sigma'(v)})^{T}, \quad r_2(v)=(1,-\sqrt{\sigma'(v)})^{T}.
\end{aligned}
\end{equation}
For the Euler system, the condition $v>0$ ensures that system is strictly hyperbolic. \eqref{sg} implies both two characteristics are genuine nonlinear, thus, only single wave pattern will form in a characteristic field and only shock wave and rarefaction wave will form \cite{JS}. However, when it comes to the viscoelastic system, both the 1-characteristic field and the 2-characteristic field are neither genuine nonlinear nor linearly degenerate. So that the multi-wave patterns will be generated on the same characteristic field. For simplicity, we analyze the multi-wave patterns on the 2-characteristic field. That is, there exists an intermediate state $v_{\ast}\in[v_-,v_+]$ such that the constant states $(v_{\pm},u_{\pm})$ as the far field states satisfy the following RH(Rankine-Hugoniot) condition:
\begin{equation}\label{RH}
\left\{
\begin{aligned}
&-\lambda_{2}(v_{\ast})(v_+-v_-)-(u_+-u_-)=0,\\
&-\lambda_{2}(v_{\ast})(u_+-u_-)-(\sigma(v_+)-\sigma(v_-))=0, v_{\ast}\in[v_-,v_+].
\end{aligned}
\right.
\end{equation}
Recalling the stress function $\sigma(v)$ we concerned here satisfies the following picture,
\begin{figure}[h]
    \centering
    \includegraphics[width=0.38\linewidth]{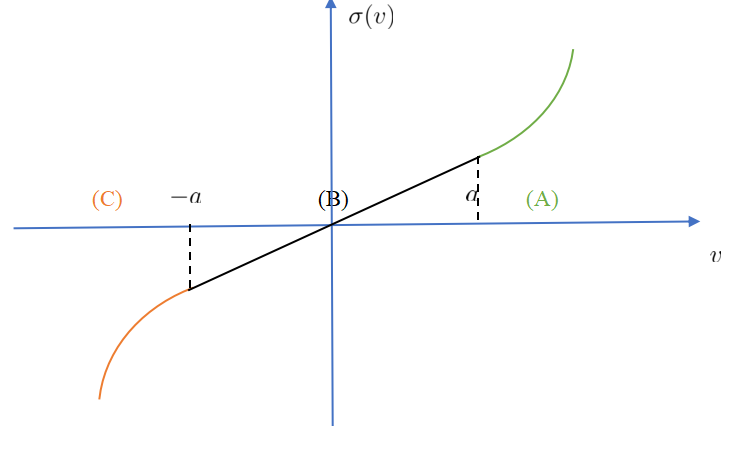}
\end{figure}
where $v_{-}$ and $v_+$ are located in two different intervals among three regions $(A),\ (B)$ and $(C)$. Thus multiple wave patterns will form in the same characteristic field. Specifically, there are three cases:  

\begin{equation}\label{cases}
\begin{aligned}
&Case\ 1:\ v_-\in(B),v_+\in(A): \text{a contact wave combines with a rarefaction wave}.\\
&Case\ 2:\ v_-\in(C),v_+\in(B):\text{a shock wave}.\\
&Case\ 3:\ v_-\in(C),v_+\in(A):\text{a shock wave or the combination of shock wave and rarefaction wave, }\\
&\quad\quad\quad\quad   \text{which depends on the position of $v_+$}.
\end{aligned}
\end{equation}
We are concerned about the first case in this paper, and the other two cases will be studied in the forthcoming paper. In the next subsection, we will construct the approximate profile for case 1.

\subsection{Construction of multi-wave consisting of rarefaction wave and contact wave}\label{Sec2.2}

In this subsection, we will give the explicit form of rarefaction wave and contact wave. Note that these two wave patterns overlap with each other in the same characteristic field, the interaction of which is very complex and we need more information about them.

\subsubsection{The approximate rarefaction wave.}

Firstly, we give the explicit form of the rarefaction wave. According to Matsumura-Nishihara \cite{MN1}, the 2-rarefaction wave $U^r=(v^r,u^r)$ could be constructed as follows: 

Recall the characteristic value $\sqrt{\sigma'(v^r)}$ satisfies the nonviscous Burgers equation with the far field states $w_{\pm}\in\mathbb{R}, (w_-<w_+)$:

\begin{align}\label{eqW}
\left\{
\begin{aligned}
&w_t+ww_x=0,\quad t>0,x\in\mathbb{R}.\\
&w(x,0)
:=\begin{cases}
w_-,\quad\quad (x<0),\\
w_+,\quad\quad (x>0).
\end{cases}
\end{aligned}
\right.
\end{align}
We denote the weak solution of \eqref{eqW} by $w^r(\frac{x}{t};w_-,w_+)$, which is 

\begin{equation}\label{wr}
w^r(\frac{x}{t};w_-,w_+):=
\begin{cases}
w_-,\quad (x\leq w_-t),\\
\frac{x}{t},\quad (w_-t\leq x\leq w_+t),\\
w_+,\quad (x\geq w_+t).
\end{cases}
\end{equation}
Order the approximate solution $w(x,t;w_-,w_+)$ which connect the far field states $w_{\pm}$,
\begin{align}\label{eqda}
\left\{
\begin{aligned}
&w_t+ww_x=0,\quad t>0,x\in\mathbb{R},\\
&w(x,0)=w_0(x)=\frac{w_{+}+w_{-}}{2}+\frac{w_{+}-w_{-}}{2} \tanh x,
\end{aligned}
\right.
\end{align}
where $w_-=\lambda_2(a),w_+=\lambda_2(v_+),$ and $a$ is from \eqref{sigma-1}.

Our rarefaction wave $U^r=(v^r,u^r)$ is defined as:
\begin{equation}\label{ur}
\begin{aligned}
&v^r(x,t;a,v_+):=(\lambda_2)^{-1}(w(x,t;\lambda_2(a),\lambda_2(v_+)),\\
&u^r(x,t;u_a,u_+):=u_{a}-\int_{a}^{v^r}\lambda_2(s)ds.
\end{aligned}   
\end{equation}
We could determine that

\begin{equation}
\begin{aligned}
u_a=u_++\int_{a}^{v_+}\lambda_2(s)ds.
\end{aligned}
\end{equation}

Then the rarefaction wave $U^r=(v^r,u^r)$ \eqref{ur} satisfies the following Cauchy problem:

\begin{equation}\label{Cauchy}
\begin{cases}
&v^{r}_t-u^{r}_x=0, \quad  t>0,x\in\mathbb{R},\\
&u^{r}_t-(\sigma(v^r))_x=0,\\
&(v^r,u^r)(x,0)=(v^r_0,u^r_0)(x)=
\begin{cases}
(a,u_a),\quad x\rightarrow -\infty,\\
(v_+,u_+),\quad x\rightarrow +\infty.
\end{cases}
\end{cases}
\end{equation}

\begin{Lem}[Properties of $U^r$]Assume $a<v_+,$ and $\sigma(v)\in C^3[a,+\infty),\sigma''(v)>0$ on $[a,v_+]$. $U^r$ defined by \eqref{ur} is the unique $C^2$ solution of the Cauchy problem \eqref{Cauchy} and it has following properties.

(1) $a<v^r(x,t)<v_+, u_+<u^r<u_{a},\partial_x v^r(x,t)>0, \partial_{x}u^r(x,t)<0,~u^r_x=O(1)v^r_x.$

(2) For $1\leq p\leq +\infty,$ there exists a positive constant $C_p$ such that 

\begin{equation}\label{Ur2}
\begin{aligned}
&\|\partial_{x}U^r(t)\|_{L^p}\leq C_p(1+t)^{-1+\frac{1}{p}},\quad (\forall t\geq 0)\\
&\|\partial_{x}^2U^r(t)\|_{L^p}\leq C_p(1+t)^{-1},\quad (\forall t\geq 0).
\end{aligned}
\end{equation}

(3) For $\forall\epsilon\in(0,1),$ there exists a positive constant $C_{\epsilon}$ such that 

\begin{equation}\label{Ur3}
\begin{aligned}
&|w(x,t)-\lambda_2(v_+)|\leq C_{\epsilon}(1+t)^{-1+\epsilon}e^{-\epsilon|x-\sqrt{\sigma'(v_+)}t|},\quad (t\geq 0, x\geq \sqrt{\sigma'(v_+)}t)\\
&|w(x,t)-\lambda_2(a)|\leq C_{\epsilon}(1+t)^{-1+\epsilon}e^{-\epsilon|x-\sqrt{\sigma'(a)}t|},\quad (t\geq 0, x\leq \sqrt{\sigma'(a)}t)\\
&|w^r(\frac{x}{t})-w(x,t)|\leq C_{\epsilon}(1+t)^{-1+\epsilon},\quad (t\geq 1, \sqrt{\sigma'(a)}t\leq x\leq \sqrt{\sigma'(v_+)}t).
\end{aligned}
\end{equation}

\end{Lem}

\subsubsection{The  contact wave.}

Instead of the nonlinear diffusion wave as that in \cite{HLM,HMX}, in this part, we construct a set of linear diffusion waves to approximate the contact wave so that we can obtain the explicit expressions of contact wave rather than its asymptotic state.

Note that the viscoelastic system \eqref{eq} is a hyperbolic-parabolic coupled system, the construction of the contact wave is not simple since diffusion waves are the solutions of parabolic equations. We need to achieve cancellations so that the errors caused by viscosity are under control. To do this, we first diagonalize system \eqref{eq} when $v\in[-a,a]$, then $\sqrt{\sigma'(v)}\equiv \sqrt{b}$. Then we define the new functions $(w_1,w_2)(x,t)$ which satisfying

\begin{equation}\label{v,u}
\left\{
\begin{aligned}
&w_1+w_2=v,\\
&\sqrt{\sigma'(v)}(w_1-w_2)=u.
\end{aligned}
\right.
\end{equation}
According to \eqref{eq}, we have

\begin{equation}\label{w1w2}
\left\{
\begin{aligned}
&w_{1t}-\sqrt{\sigma'(v)}w_{1x}=\frac{\mu}{2}(w_1-w_2)_{xx},\\
&w_{2t}+\sqrt{\sigma'(v)}w_{2x}=\frac{\mu}{2}(w_2-w_1)_{xx}.
\end{aligned}
\right.
\end{equation}
In order to make the functions $w_{1}$ and $w_{2}$ satisfy \eqref{w1w2} as much as possible, we give the forms of $w_{1}$ and $w_{2}$ as follows

\begin{equation}\label{Xi}
\left\{
\begin{aligned}
&w_1=\Xi_1+\frac{\mu}{4\sqrt{\sigma'(v)}}\Xi_{2x},\\
&w_2=\Xi_2-\frac{\mu}{4\sqrt{\sigma'(v)}}\Xi_{1x}.\\
\end{aligned}
\right.
\end{equation}
where $\Xi_1,\Xi_2$ satisfy the following linear parabolic equations, $\sqrt{\sigma'(v)}=\sqrt{b},$

\begin{equation}\label{Xi12}
\left\{
\begin{aligned}
&\Xi_{1t}-\sqrt{\sigma'(v)}\Xi_{1x}=\frac{\mu}{2}\Xi_{1xx},\quad  \Xi_1(-\infty,t)=0,\Xi_1(+\infty,t)=0,\\
&\Xi_{2t}+\sqrt{\sigma'(v)}\Xi_{2x}=\frac{\mu}{2}\Xi_{2xx},\quad \Xi_2(-\infty,t)=v_-,\Xi_2(+\infty,t)=a.
\end{aligned}
\right.
\end{equation}
From $RH$ condition, we know,

\begin{equation}
\begin{aligned}
-\sqrt{b}(a-v_-)-(u_{a}-u_-)=0,\\
-\sqrt{b}(u_{a}-u_-)-(\sigma(a)-\sigma(v_-))=0.
\end{aligned}
\end{equation}
Then without loss of generality, we set 
\begin{align}
 \frac{u_-+\sqrt{b}v_-}{2\sqrt{b}}=\frac{u_{a}+\sqrt{b}a}{2\sqrt{b}}=0, \quad \frac{\sqrt{b}v_--u_-}{2\sqrt{b}}=v_-,\frac{\sqrt{b}a-u_{a}}{2\sqrt{b}}=a.
\end{align}
Then it follows that,
\begin{equation}
\begin{aligned}
&\Xi_{1\pm}=\Xi_1(\pm\infty,t)=w_{1\pm}=w_1(\pm\infty,t)=0,\\
&\Xi_{2-}=\Xi_2(-\infty,t)=w_{2_-}=w_2(-\infty,t)=v_-, \\
&\Xi_{2+}=\Xi_2(+\infty,t)=w_{2+}=w_2(+\infty,t)=a.
\end{aligned}
\end{equation}
Solving \eqref{Xi12}, the explicit expression of $\Xi_1,\Xi_2$ are 
\begin{equation}\label{XiXi}
\begin{aligned}
&\Xi_1=0,\\
&\Xi_2=v_-+\frac{a-v_-}{\sqrt{\pi}}\int^{\frac{x-\sqrt{b}t}{\sqrt{2\mu (1+t)}}}_{-\infty}e^{-\xi^2}d\xi.
\end{aligned}
\end{equation}
We could get the viscous contact wave $U^c=(v^c(\frac{x-\sqrt{b}t}{\sqrt{1+t}};v_-,a),u^c(\frac{x-\sqrt{b}t}{\sqrt{1+t}};u_-,u_a))$ as 

\begin{equation}\label{uc}
v^c=\Xi_{2}+\frac{\mu}{4\sqrt{b}}\Xi_{2x},\quad u^c=\sqrt{b}(\frac{\mu}{4\sqrt{b}}\Xi_{2x}-\Xi_2).
\end{equation}
Then we could derive the equation of viscous contact wave  ${U}^{c}=(v^c,u^c)$ as follows:
\begin{equation}\label{uceq}
\left\{
\begin{aligned}
&v^c_t-u^c_x=Q_1,\\
&u^c_t-(\sigma(v^c))_x=\mu u^c_{xx}+Q_2,
\end{aligned}
\right.
\end{equation}
where the remainder terms are

\begin{equation}
\begin{aligned}
&Q_1=\frac{\mu}{2}\frac{\mu}{4\sqrt{b}}\Xi_{2xxx},\quad Q_2=-\sqrt{b}\frac{\mu}{2}\frac{\mu}{4\sqrt{b}}\Xi_{2xxx},\\
&|(Q_1,Q_2)|\lesssim O(\delta)(1+t)^{-\frac{3}{2}}e^{-\frac{c(x-\sqrt{b}t)^2}{2\mu(1+t)}}.
\end{aligned}
\end{equation}
Now we define the composite wave of rarefaction wave $U^r=(v^r,u^r)$ and viscous contact wave $U^c=(v^c,u^c)$ as $\hat{U}=(\hat{v},\hat{u})$, 

\begin{equation}\label{hatu}
\begin{aligned}
&\hat{v}(x,t)=v^c(\frac{x-\sqrt{b}t}{\sqrt{1+t}};v_-,a)+v^r(\frac{x}{1+t};a,v_+)-a,\\
&\hat{u}(x,t)=u^c(\frac{x-\sqrt{b}t}{\sqrt{1+t}};u_-,u_{a})+u^r(\frac{x}{1+t};u_{a},u_+)-u_{a}.
\end{aligned}
\end{equation}
By \eqref{Cauchy} and \eqref{uceq},  we derive the equation of $\hat{U}$,

\begin{equation}\label{hatueq}
\left\{
\begin{aligned}
&\hat{v}_t-\hat{u}_x=Q_1,\\
&\hat{u}_t-(\sigma(\hat{v}))_x=\mu\hat{u}_{xx}+H.
\end{aligned}
\right.
\end{equation}
where $H$ is defined as 

\begin{equation}\label{H}
\begin{aligned}
H=-(\sigma(\hat{v})-\sigma(v^c)-\sigma(v^r))_x-\mu u^r_{xx}+Q_2.
\end{aligned}
\end{equation}
The initial data of \eqref{hatueq} is 

\begin{equation}\label{hatin}
\begin{aligned}
(\hat{v},\hat{u})(x,0)=(\hat{v}_0,\hat{u}_0)(x)\rightarrow (v_{\pm},u_{\pm}),\quad \textbf{as}\quad x\rightarrow\pm\infty.
\end{aligned}
\end{equation}
Finally we are able to give the mathematical description of our main theorem.
\subsection{Main results}
We give the main theorem in this paper.

\begin{Thm}\label{mt}Assume that function $\sigma(v)$ satisfies \eqref{sigma-1}, the far field states $(v_{\pm},u_{\pm})$ satisfies the $RH$ condition \eqref{RH}. Moreover, they belong to Case 1. If there exists a constant $\epsilon_0>0$ such that if
\begin{equation}\label{t1}
\begin{aligned}
 \delta+\|(v_0,u_0)-(\hat{v}_0,\hat{u}_0)\|_{H^1(\mathbb{R})}\leq \epsilon_0, 
\end{aligned}
\end{equation}
where $\delta=|v_+-v_-|$, then there exists a unique global smooth solution $(v,u)(x,t)$ of \eqref{eq},\eqref{inta} satisfying
	\begin{equation}\label{st}
	\begin{aligned}
\lim_{t\rightarrow\infty}\sup_{x\in\mathbb{R}}|(v-\hat{v},u-\hat{u})(x,t)|=0.
	\end{aligned}
	\end{equation}
\end{Thm}
\begin{Rem}
    To the best of our knowledge, \cref{mt} is the first stability result of composite wave patterns generated in the same characteristic field of systems not limited to the system of viscoelastic.
\end{Rem}
\begin{Rem}
    From \cref{hatu}, one knows the two different wave patterns will not separate from each other even if the time goes into infinity. This is essentially different from the previous remarkable results of composite wave patterns such as \cite{HLM,HM,KVW1,KVW2}. The interaction is much more complicated and is the key to this paper. We analyze the interaction surface carefully in \cref{sec3}.
\end{Rem}
\begin{Rem}
    The method to construct the viscous contact wave is completely different from \cite{HMX,HXY}, in which the nonlinear diffusion wave is used to approximate the viscous contact wave in a subtle way. The contact wave in our paper consists of a series of linear diffusion waves thus the explicit expressions of contact wave are available. Note that the construction in this paper is based on the diagonalized system and we claim that this method can be applied to more general hyperbolic-parabolic systems.
\end{Rem}
\begin{Rem}
    The leading term for viscous contact wave is linear diffusion wave, thus the critical terms 
    \begin{align}
        \int_0^T\int_{\R}(1+\tau)^{-1}e^{-\frac{c\eta^2}{1+\tau}}(\phi^2,\psi^2) d\eta d\tau
    \end{align}
    will appear in the stability analysis. Motivated by the estimates of heat kernel introduced by \cite{HLM}, we developed a new Poincar\'e type inequality to close the fundamental energy inequality. Note that the flux is generally non-convex, the estimate is much more complex and this method can be applied to other systems with non-convex flux, see \cref{4.2} for more details.
\end{Rem}

\section{The interaction of the contact wave and the rarefaction wave}\label{sec3}

When $-a<v_-<a<v_+,$ we know that for $\forall t\geq 0,$ there exists $X_1(t)\in\mathbb{R}$ as the curve connecting two wave patterns based on the monotonicity of functions $\Xi_2,v^r$.

\begin{equation}
(\hat{v}-\frac{\mu}{4\sqrt{b}}\Xi_{2x})(X_1(t),t)=\Xi_2(X_1(t),t)+v^{r}(X_1(t),t)-a=a\quad (t\geq 0)
\end{equation}
From \eqref{XiXi} it yields that

\begin{equation}\label{X_1(t)}
v^{r}(X_1(t),t)-a=a-\Xi_2(X_1(t),t)=\frac{a-v_-}{\sqrt{\pi}}\int_{\frac{X_1(t)-\sqrt{b}t}{\sqrt{2\mu (1+t)}}}^{+\infty}e^{-\xi^2}d\xi\quad (t\geq 0)
\end{equation}
Similar as \cite{MY}, $X_1(t)$ has following properties,

\begin{Lem}(\cite{MY})\label{xl}
Assume $-a<v_-<a<v_+$, the following properties hold

(1) There exists a positive $T_0$, such that, for any $t>T_0$,
	\begin{align}\label{l1}
	\sqrt{b}(1+t)+\sqrt{2\mu(1+t)}\leq X_1(t)\leq \lambda_{2}(v_+)(1+t).
	\end{align}

 (2) Define $X_1(t)$ as in (\ref{X_1(t)}), for $t>T_0$, we have

	\begin{align}\label{l2}
	\left|(\lambda_2)^{-1}\left(\frac{X_1(t)}{1+t}\right)-a-\frac{a-v_-}{\sqrt{\pi}} \int_{\frac{X_1(t)-\sqrt{b}t}{\sqrt{2\mu (1+t)}}}^{\infty} \mathrm{e}^{-\xi^{2}} \mathrm{d} \xi\right| \leq C(1+t)^{-\frac{3}{4}}.
	\end{align}

 (3) For any positive constant $\epsilon\in(0,1)$, there exists a constant $C_\epsilon>0,t\geq T_0$ such that
	
\begin{equation}\label{l3}
    \begin{aligned}
	 &X_1(t)\leq \sqrt{b}(1+t)+\left(C+ln(1+t)^{\frac{1}{2}}\right)^{\frac{1}{2}}\sqrt{2\mu(1+t)},\\
	 &X_1(t)\geq \sqrt{b}(1+t)+\left(C_\epsilon+ln(1+t)^{\frac{1}{2(1+\epsilon)}}\right)^{\frac{1}{2}}\sqrt{2\mu(1+t)},
   \end{aligned}
\end{equation}
where $ln\ t= \log_e (t)$.

\end{Lem}

\begin{proof}
\noindent{step 1:} Order $x_l(t)=\sqrt{b}(1+t)+\sqrt{2\mu(1+t)},$ from  \eqref{Ur3},\eqref{XiXi} as $t\rightarrow+\infty,$

\begin{equation}
\begin{aligned}
\hat{v}(x_l(t),t)=v^c(x_l(t),t)+v^r(x_l(t),t)-a\rightarrow v_-+\frac{a-v_-}{\sqrt{\pi}}\int_{-\infty}^{1}e^{-\xi^2}d\xi<a
\end{aligned}
\end{equation}

Order $x_r(t)=\lambda_{2}(v_+)(1+t),$ as $t\rightarrow+\infty,$

\begin{equation}
\begin{aligned}
\hat{v}(x_r(t),t)=v^c(x_r(t),t)+v^r(x_r(t),t)-a\rightarrow v_+>a
\end{aligned}
\end{equation}
Then we could get \eqref{l1}.

\noindent{step 2: } From \eqref{Ur3} and \eqref{X_1(t)}, because 

\begin{equation}
|w^r(\frac{X_1(t)}{1+t})-w(X_1(t),t)|=|w^r(\frac{X_1(t)}{1+t})-\sqrt{b}+\sqrt{b}-w(X_1(t),t)|\leq  C_{\epsilon}(1+t)^{-1+\epsilon},
\end{equation}
and $\lambda_2(v)$ is invertible  on $[a,v_+],$ $\lambda_2(a)=\sqrt{b},$ it yields that
\begin{equation}
\begin{aligned}
|(\lambda_2)^{-1}\left(\frac{X_1(t)}{1+t}\right)-a-\frac{a-v_-}{\sqrt{\pi}} \int_{\frac{X_1(t)-\sqrt{b}t}{\sqrt{2\mu (1+t)}}}^{\infty} \mathrm{e}^{-\xi^{2}} \mathrm{d} \xi|
\leq C_{\epsilon}(1+t)^{-1+\epsilon}.
\end{aligned}
\end{equation}
Take $\epsilon=\frac{1}{4},$ we could get \eqref{l2} immediately.

\noindent{step 3: } Similar as the methods in \cite{MY},
putting $Y_1(t):=\frac{X_1(t)-\sqrt{b}(1+t)}{\sqrt{2\mu(1+t)}},$ noting that $Y_1(t)\geq 1$ by \eqref{l1}. By utilizing the invertibility of $\lambda_2(v)$ on $[a,v_+]$, from \eqref{l2}, we have 

\begin{equation}
\begin{aligned}
\frac{Y_1(t)}{\sqrt{1+t}}\leq & C\int_{Y(t)}^{+\infty}e^{-\xi^2}d\xi+C(1+t)^{-\frac{3}{4}}\\
\leq &C\int_{Y(t)}^{+\infty}\xi e^{-\xi^2}d\xi+C(1+t)^{-\frac{3}{4}}\\
\leq &Ce^{-Y(t)^2}+C(1+t)^{-\frac{3}{4}},\quad (t\geq T_0)
\end{aligned}
\end{equation}
This gives that 

\begin{equation}
\begin{aligned}
Y_1(t):=\frac{X_1(t)-\sqrt{b}(1+t)}{\sqrt{2\mu(1+t)}}\leq (C+log(1+t)^{\frac{1}{2}})^{\frac{1}{2}}.
\end{aligned}
\end{equation}
So we prove the first inequality of \eqref{l3}. As to prove the another inequality of \eqref{l3}, the method is same as in \cite{MY}, we omit the details.

\end{proof}

\begin{Lem}\label{le3.2}The remainder term $H$ given by \eqref{H} has the following $L^p$-estimates in time. 

(i) For $\forall \epsilon\in (0,1),$ there exist positive constant $C_{\epsilon}$ and $T_0>0$ such that

\begin{equation}\label{HL1}
\|H\|_{L^1}\leq C_{\epsilon}(1+t)^{-\frac{1}{2}-\frac{1}{2(1+\epsilon)}}\sqrt{\ln(2+t)},\forall \epsilon\in (0,1),\quad (t\geq T_0)
\end{equation}

(ii) For $p=\infty,$ there exists a positive constant $C$ such that

\begin{equation}\label{Hinfty}
\|H\|_{L^{\infty}}\leq C(1+t)^{-\frac{1}{2}},\quad (t\geq 0)
\end{equation}

(iii) Especially from the interpolation inequality, we have 

\begin{equation}\label{lp}
\begin{aligned}
\|H\|_{L^p}\leq \|H\|_{L^1}^{\frac{1}{p}}\|H\|_{L^{\infty}}^{\frac{p-1}{p}}
\leq C_{\epsilon}(1+t)^{-\frac{1}{2}-\frac{1}{2p(1+\epsilon)}}(\ln(2+t))^{\frac{1}{2p}}.
\end{aligned}
\end{equation}

\end{Lem}

\begin{proof}
\noindent{step 1:} Obviously, the terms $\|(u^r_{xx},Q_2)\|_{L^1}\leq C(1+t)^{-1}$ satisfying the decay estimates \eqref{HL1}. From the definition of \eqref{H}, we only need to prove 
\begin{equation}\label{H1}
\begin{aligned}
&\int_{-\infty}^{+\infty}|\sigma'(\hat{v})-\sigma'(v^c)||v^c_x|dx\leq C_{\epsilon}(1+t)^{-\frac{1}{2}-\frac{1}{2(1+\epsilon)}}\sqrt{\ln(2+t)},\\
&\int_{-\infty}^{+\infty}|\sigma'(\hat{v})-\sigma'(v^r)||v^r_x|dx\leq C_{\epsilon}(1+t)^{-\frac{1}{2}-\frac{1}{2(1+\epsilon)}}\sqrt{\ln(2+t)}.
\end{aligned}
\end{equation}
To prove the first inequality, note that $v^c_x=\Xi_{2x}+\frac{\mu}{4\sqrt{b}}\Xi_{2xx},$ we divide it into two parts,

\begin{equation}
\begin{aligned}
&\int_{-\infty}^{+\infty}|\sigma'(\hat{v})-\sigma'(v^c)||v^c_x|dx\leq \int_{-\infty}^{X_1(t)}|\sigma'(\hat{v})-\sigma'(v^c)||v^c_x|dx+\int_{X_1(t)}^{+\infty}|\sigma'(\hat{v})-\sigma'(v^c)||v^c_x|dx
\end{aligned}
\end{equation}
Order $\hat{v}_1=(\hat{v}-\frac{\mu}{4\sqrt{b}}\Xi_{2x})=\Xi_2+v^r-a,$ note that when $x\leq X_1(t)$, $\hat{v}_1\leq a,$ then there are three cases, $\hat{v}\geq v^c\geq a\geq \hat{v}_1,$ or $\hat{v}\geq a\geq v^c,$ or $a\geq\hat{v}\geq v^c$,

\begin{equation}\label{iqq}
\begin{aligned}
\int_{-\infty}^{X_1(t)}|\sigma'(\hat{v})-\sigma'(v^c)||v^c_x|dx\leq &C\int_{-\infty}^{X_1(t)}\bigg\{|\sigma'(\hat{v})-\sigma'(a)|+|\sigma'(a)-\sigma'(v^c)|\bigg\}(|\Xi_{2x}|+|\Xi_{2xx}|) dx\\
\leq &\int_{-\infty}^{X_1(t)}\bigg\{|\sigma'(\hat{v})-\sigma'(a)|+|\sigma'(a)-\sigma'(v^c)|\bigg\}(|\Xi_{2x}|+|\Xi_{2xx}|) dx.
\end{aligned}
\end{equation}
Since $\sigma''(v)=0,~v<a$, when $\hat{v}\leq a$ and $v^c\leq a,$  the integral is zero, so we only consider the case $\hat{v}\geq v^c>a.$ 
	
We know that $\hat{v}_1\rightarrow v_{\pm}$ as $x\rightarrow\pm\infty$ and $\hat{v}_{1x}>0$, so there exists a unique $Z_1(t)\in\mathbb{R}$ such that 

\begin{equation}\label{ZZ}
\hat{v}_1(Z_1(t),t)=a-\frac{\frac{\mu}{4\sqrt{b}}(a-v_-)}{\sqrt{\pi}\sqrt{2\mu(1+t)}}>v_-,\quad t\geq T_0.
\end{equation}

Then for $x\leq Z_1(t),$ we have

\begin{equation}\label{Z1}
\begin{aligned}
&\Xi_{2}(x,t)\leq \hat{v}_1(x,t)\leq a-\frac{\frac{\mu}{4\sqrt{b}}(a-v_-)}{\sqrt{\pi}\sqrt{2\mu(1+t)}}\leq a-\frac{\frac{\mu}{4\sqrt{b}}(a-v_-)}{\sqrt{\pi}\sqrt{2\mu(1+t)}}e^{-\frac{(x-\sqrt{b}t)^2}{2\mu(1+t)}}= a-\frac{\mu}{4\sqrt{b}}\Xi_{2x},\\
&v^c(x,t)\leq \hat{v}(x,t)=\hat{v}_1+\frac{\mu}{4\sqrt{b}}\Xi_{2x}\leq a, \quad x\leq Z_1(t),
\end{aligned}
\end{equation}
then we immediately get that 

\begin{equation}\label{Z1-1}
\begin{aligned}
0&\leq v^{r}(X_1(t),t)-a=a-\Xi_2(X_1(t),t)=\frac{a-v_-}{\sqrt{\pi}}\int_{\frac{X_1(t)-\sqrt{b}t}{\sqrt{2\mu (1+t)}}}^{+\infty}e^{-\xi^2}d\xi,\\
0&\leq v^{r}(Z_1(t),t)-a\\&=\hat{v}(Z_1(t),t)-\Xi_2(Z_1(t),t)-\frac{\mu}{4\sqrt{b}}\Xi_{2x}(Z_1(t),t)\\
&\leq a-\Xi_2(Z_1(t),t)-\frac{\mu}{4\sqrt{b}}\Xi_{2x}(Z_1(t),t)\\&=\frac{a-v_-}{\sqrt{\pi}}\int_{\frac{Z_1(t)-\sqrt{b}t}{\sqrt{2\mu (1+t)}}}^{+\infty}e^{-\xi^2}d\xi-\frac{\frac{\mu}{4\sqrt{b}}(a-v_-)}{\sqrt{\pi}\sqrt{2\mu(1+t)}}e^{-\frac{(Z_1(t)-\sqrt{b}t)^2}{2\mu(1+t)}}.
\end{aligned}
\end{equation}
Obviously, $Z_1(t)\leq X_1(t), v^{r}(Z_1(t),t)\leq v^{r}(X_1(t),t),$ this gives that

\begin{equation}
\begin{aligned}
&0\leq \frac{a-v_-}{\sqrt{\pi}}\int_{\frac{Z_1(t)-\sqrt{b}t}{\sqrt{2\mu (1+t)}}}^{+\infty}e^{-\xi^2}d\xi-\frac{\frac{\mu}{4\sqrt{b}}(a-v_-)}{\sqrt{\pi}\sqrt{2\mu(1+t)}}e^{-\frac{(Z_1(t)-\sqrt{b}t)^2}{2\mu(1+t)}}\leq \frac{a-v_-}{\sqrt{\pi}}\int_{\frac{X_1(t)-\sqrt{b}t}{\sqrt{2\mu (1+t)}}}^{+\infty}e^{-\xi^2}d\xi,
\end{aligned}
\end{equation}
thus
\begin{equation}\label{Z1-2}
\begin{aligned}
&\frac{a-v_-}{\sqrt{\pi}}\int_{Z_1(t)}^{X_1(t)}e^{-\frac{(x-\sqrt{b}t)^2}{2\mu(1+t)}}dx=\frac{a-v_-}{\sqrt{\pi}}\int_{\frac{Z_1(t)-\sqrt{b}t}{\sqrt{2\mu (1+t)}}}^{\frac{X_1(t)-\sqrt{b}t}{\sqrt{2\mu (1+t)}}}e^{-\xi^2}\sqrt{2\mu(1+t)}d\xi\leq \frac{\mu(a-v_-)}{4\sqrt{b\pi}}e^{-\frac{(Z_1(t)-\sqrt{b}t)^2}{2\mu(1+t)}}.
\end{aligned}
\end{equation}

Using \eqref{Z1},\eqref{Z1-2}, and note that $\hat{v}-a\leq \hat{v}-\hat{v}_1=\frac{\mu}{4\sqrt{b}}\Xi_{2x}, x\leq X_1(t)$, it yields that
\begin{equation}\label{Z1-3}
\begin{aligned}
&\int_{-\infty}^{X_1(t)}\bigg\{|\sigma'(\hat{v})-\sigma'(a)|\bigg\}(|\Xi_{2x}|+|\Xi_{2xx}|) dx\\
=&\int_{Z_1(t)}^{X_1(t)}\bigg(\sigma''(\hat{v})(\hat{v}-a)+O(\hat{v}-a)^2\bigg)(|\Xi_{2x}|+|\Xi_{2xx}|) dx\\
\leq & C\int_{Z_1(t)}^{X_1(t)}|\Xi_{2x}|^2+|\Xi_{2x}||\Xi_{2xx}|dx\\
\leq &C(1+t)^{-1}\int_{Z_1(t)}^{X_1(t)}e^{-\frac{(x-\sqrt{b}t)^2}{2\mu(1+t)}}dx+C(1+t)^{-\frac{3}{2}}\int_{-\infty}^{+\infty}e^{-\frac{c(x-\sqrt{b}t)^2}{(1+t)}}dx\\
\leq &C(1+t)^{-1}.
\end{aligned}  
\end{equation}
Similarly as \eqref{Z1-3}, from \eqref{Z1}, when $v^c>a,$ we have

\begin{equation}\label{Z1-4}
\begin{aligned}
\int_{-\infty}^{X_1(t)}\bigg\{|\sigma'(a)-\sigma'(v^c)|\bigg\}(|\Xi_{2x}|+|\Xi_{2xx}|) dx
\leq C(1+t)^{-1}.
\end{aligned}  
\end{equation}
As for $x\geq X_1(t),$ we have $\hat{v}\geq \hat{v}_1\geq v^c\geq a,$ $\hat{v}\geq v^c\geq \hat{v}_1\geq a,$ or $\hat{v}\geq \hat{v}_1\geq a\geq v^c,$ and by \eqref{l3},

\begin{equation}\label{Xi2x}
	\Xi_{2xx}=O(\delta)(1+t)^{-1}\frac{-(x-\sqrt{b}t)}{\sqrt{2\mu(1+t)}}e^{-\frac{(x-\sqrt{b}t)^2}{2\mu(1+t)}}<0,
\end{equation}
thus
\begin{equation}\label{iq1}
\begin{aligned}
\int_{X_1(t)}^{+\infty}|\sigma'(\hat{v})-\sigma'(v^c)||v^c_x|dx\leq &C\int_{X_1(t)}^{+\infty}(\sigma'(\hat{v})-\sigma'(v^c))(\Xi_{2x}-\Xi_{2xx})dx\\
\leq &C\int_{X_1(t)}^{+\infty}\big\{|\hat{v}-\hat{v}_1|+|\hat{v}_1-v^c|\big\}(\Xi_{2x}-\Xi_{2xx}) dx\\
\leq &C\int_{X_1(t)}^{+\infty}\big\{|\Xi_{2x}|+(v^r-a)\big\}(\Xi_{2x}-\Xi_{2xx}) dx\\
\leq &C\int_{X_1(t)}^{+\infty}|\Xi_{2x}|^2+|\Xi_{2xx}||\Xi_{2x}|dx-C\int_{X_1(t)}^{+\infty}v^r_x(\Xi_2-a-\Xi_{2x})dx\\&~~+(v^r-a)(\Xi_2-a-\Xi_{2x})\mid_{X_1(t)}^{+\infty}.
\end{aligned}
\end{equation}
Here 

\begin{equation}\label{iq2}
\begin{aligned}
\int_{X_1(t)}^{+\infty}|\Xi_{2x}|^2+|\Xi_{2xx}||\Xi_{2x}|dx
\leq \|(\Xi_{2x},\Xi_{2xx})\|_{L^{\infty}}\int_{X_1(t)}^{+\infty}\Xi_{2x}dx
\leq  C(1+t)^{-\frac{1}{2}}|\Xi_2(X_1(t),t)-a|.
\end{aligned}
\end{equation}
By \eqref{l2} and \eqref{l3}, we have
\begin{equation}\label{vraxia}
	\begin{aligned}
		|v^r(X_1(t),t)-a|=|\Xi_2(X_1(t),t)-a|\leq C\bigg|\frac{X_{1}(t)-\sqrt{b}(1+t)}{1+t}\bigg|+C(1+t)^{-\frac{3}{4}}\leq C\bigg(\frac{\ln(2+t)}{1+t}\bigg)^{\frac{1}{2}},
	\end{aligned}
\end{equation}
thus
\begin{equation}\label{iq3}
	\begin{aligned}
		|(v^r(X_1(t),t)-a)||(\Xi_2(X_1(t),t)-a)|\leq C(1+t)^{-1}\ln(2+t).
	\end{aligned}
\end{equation}
For the remaining term in \eqref{iq1},
\begin{equation}\label{iq4}
\begin{aligned}
\int_{X_1(t)}^{+\infty}v^r_x(a-\Xi_2+\Xi_{2x})dx
\leq \|v^r_x\|_{L^{\infty}}\int_{X_1(t)}^{+\infty}(a-\Xi_2+\Xi_{2x})dx
\leq C_{\epsilon}(1+t)^{-\frac{1}{2}-\frac{1}{2(1+\epsilon)}},
\end{aligned}
\end{equation}
where

\begin{equation}\label{iq5}
\begin{aligned}
&\|v^r_x\|_{L^{\infty}}\int_{X_1(t)}^{+\infty}(a-\Xi_2+\Xi_{2x})dx\\
\leq &C(1+t)^{-1}\bigg\{\int_{X_1(t)}^{+\infty}\int_{\frac{x-\sqrt{b}t}{\sqrt{2\mu(1+t)}}}^{+\infty}e^{-\xi^2}d\xi dx+\int_{X_1(t)}^{+\infty}(1+t)^{-\frac{1}{2}}e^{-\frac{(x-\sqrt{b}t)^2}{2\mu(1+t)}}dx\bigg\}\\
\leq &C(1+t)^{-1}\int_{\frac{X_1(t)-\sqrt{b}t}{\sqrt{2\mu(1+t)}}}^{+\infty}e^{-\xi^2}(\xi\sqrt{2\mu(1+t)}+\sqrt{b}t-X_1(t))d\xi\\
\leq &(1+t)^{-\frac{1}{2}}e^{-\frac{(X_1(t)-\sqrt{b}t)^2}{2\mu(1+t)}}\\
\leq &C_{\epsilon}(1+t)^{-\frac{1}{2}-\frac{1}{2(1+\epsilon)}}.
\end{aligned}
\end{equation}
by using \eqref{l3}. Inserting \eqref{iq2}-\eqref{iq5} into \eqref{iq1}, we could obtain

\begin{equation}\label{iq6}
\begin{aligned}
\int_{X_1(t)}^{+\infty}|\sigma'(\hat{v})-\sigma'(v^c)||v^c_x|dx\leq C_{\epsilon}(1+t)^{-\frac{1}{2}-\frac{1}{2(1+\epsilon)}}.
\end{aligned}
\end{equation}
Combining \eqref{iqq}, \eqref{Z1-3}-\eqref{Z1-4} and \eqref{iq6}, we finish the proof of the first inequality in \eqref{H1}.

To prove the second inequality in \eqref{H1}, note that $v^r_x>0, \hat{v}(x,t)\geq a, x\geq X_1(t)$ still holds, then 

\begin{equation}\label{iq7}
\begin{aligned}
\int_{-\infty}^{+\infty}|\sigma'(\hat{v})-\sigma'(v^r)||v^r_x|dx
\leq & \int_{-\infty}^{X_1(t)}|\sigma'(\hat{v})-\sigma'(v^r)|v^r_xdx+\int_{X_1(t)}^{+\infty}|\sigma'(\hat{v})-\sigma'(v^r)|v^r_xdx.
\end{aligned}
\end{equation}

Remember that $\hat{v}_1=(\hat{v}-\frac{\mu}{4\sqrt{b}}\Xi_{2x}),$ then $\hat{v}\geq \hat{v}_1,$ 
note that $\hat{v}_1(x,t)\leq a, x\leq X_1(t),$ if $\hat{v}-a< 0,$ then $|\sigma'(\hat{v})-\sigma'(a)|=0$. 
 When $\hat{v}-a\geq 0$ and $x\leq X_1(t),$
 we have $ |\hat{v}-a|\leq |\hat{v}-\hat{v}_1|=\frac{\mu}{4\sqrt{b}}\Xi_{2x}$, thus
\begin{equation}\label{iq8}
\begin{aligned}
\int_{-\infty}^{X_1(t)}|\sigma'(\hat{v})-\sigma'(v^r)|v^r_xdx\leq &\int_{-\infty}^{X_1(t)}\big(|\sigma'(\hat{v})-\sigma'(a)|+|\sigma'(a)-\sigma'(v^r)|\big)v_x^rdx\\
\leq &C\int_{-\infty}^{X_1(t)}\big(|\hat{v}-a|+|a-v^r|\big)v^r_xdx,\\
\leq &C\int_{-\infty}^{X_1(t)}((\Xi_{2}-v_-)v^r_x)_x-(\Xi_2-v_-)v^r_{xx}+\frac{1}{2}((v^r-a)^2)_xdx\\
\leq &C\bigg\{v^r_x(X_1(t),t)+(v^r(X_1(t),t)-a)^2+\|v^r_{xx}\|_{L^1(\mathbb{R})}\bigg\}\\
\leq &C(1+t)^{-1}+C(1+t)^{-1}\ln (2+t)\\
\leq &
C(1+t)^{-1}\ln (2+t),
\end{aligned}
\end{equation}
where we have used \eqref{vraxia}. By \eqref{l3},

\begin{equation}\label{iq9}
\begin{aligned}
\int_{X_1(t)}^{+\infty}|\sigma'(\hat{v})-\sigma'(v^r)|v^r_xdx\leq &C\int_{X_1(t)}^{+\infty}|\Xi_2+\frac{\mu}{4\sqrt{b}}\Xi_{2x}-a|v_x^rdx\\
\leq &\|v_x^r\|_{L^{\infty}}\int_{X_1(t)}^{+\infty}|\Xi_2-a|+|\frac{\mu}{4\sqrt{b}}\Xi_{2x}|dx\\
\end{aligned}
\end{equation}
The estimation of \eqref{iq9} is similar as \eqref{iq5}. Then we prove \eqref{H1}.

\noindent{step 2:} Obviously, from \eqref{H}, 

\begin{equation}
\|H\|_{\infty}\leq C\|(Q_2,u^r_{xx},v^c_x,v^r_x)\|_{\infty}\leq C(1+t)^{-\frac{1}{2}}
\end{equation}

\noindent{step 3:} The interpolation inequality \eqref{lp} is standard, we omit the details.

\end{proof}

\begin{Lem}\label{le3.3}The location of $Z_1(t)$ satisfies following inequality:

\begin{equation}\label{lem3.3}
\begin{aligned}
\sqrt{b}t+(\ln(1+t)^{\beta})^{\frac{1}{2}}\sqrt{2\mu(1+t)}\leq Z_1(t)\leq X_1(t),\quad t\geq T_0
\end{aligned}
\end{equation}
for some small constant $\beta>0$, where $Z_1(t)$ is defined in \eqref{ZZ}.
 \end{Lem}

\begin{proof} Since $Z_1(t)$ satisfies  
\begin{equation}
\hat{v}_1(Z_1(t),t)=a-\frac{\frac{\mu}{4\sqrt{b}}(a-v_-)}{\sqrt{\pi}\sqrt{2\mu(1+t)}}>v_-,\quad t\geq T_0.
\end{equation}
Order $M_1(t)=\sqrt{b}t+(\ln(1+t)^{\beta})^{\frac{1}{2}}\sqrt{2\mu(1+t)},$ then 

\begin{equation}
\begin{aligned}
&\hat{v}_1(M_1(t),t)=\bigg\{\Xi_2+v^r-a\bigg\}(M_1(t),t)\\
=&v_-+\frac{a-v_-}{\sqrt{\pi}}\int_{-\infty}^{\frac{M_1(t)-\sqrt{b}t}{\sqrt{2\mu(1+t)}}}e^{-\xi^2}d\xi+v^r(M_1(t),t)-a\\
=&v_-+\frac{a-v_-}{\sqrt{\pi}}\int_{-\infty}^{(\ln(1+t)^{\beta})^{\frac{1}{2}}}e^{-\xi^2}d\xi+v^r(M_1(t),t)-a
\end{aligned}
\end{equation}
Now we show that 
$\hat{v}_1(M_1(t),t)\leq \hat{v}_1(Z_1(t),t),$ that is 
\begin{equation}\label{h1z1}
\begin{aligned}
v_-+\frac{a-v_-}{\sqrt{\pi}}\int_{-\infty}^{(\ln(1+t)^{\beta})^{\frac{1}{2}}}e^{-\xi^2}d\xi+v^r(M_1(t),t)-a\leq a-\frac{\frac{\mu}{4\sqrt{b}}(a-v_-)}{\sqrt{\pi}\sqrt{2\mu(1+t)}}
\end{aligned}
\end{equation}
for some small constant $\beta$. It is equivanlent to show 

\begin{equation}\label{in1}
\begin{aligned}
\frac{\frac{\mu}{4\sqrt{b}}(a-v_-)}{\sqrt{\pi}\sqrt{2\mu(1+t)}}\leq \frac{a-v_-}{\sqrt{\pi}}\int_{(\ln(1+t)^{\beta})^{\frac{1}{2}}}^{+\infty}e^{-\xi^2}d\xi+a-v^r(M_1(t),t)
\end{aligned}
\end{equation}
By using $\eqref{Ur3}_3$, when $t\geq T_0,$ 

\begin{equation}\label{loc1}
\begin{aligned}
&|a-v^r(M_1(t),t)|
=|\lambda^{-1}(\sqrt{b})-\lambda^{-1}(w(M_1(t),t))|\\
=&O(1)|\sqrt{b}-w(M_1(t),t)|,\quad (w^{r}(\frac{M_1(t)}{t})=\frac{M_1(t)}{t})\\
\leq &O(1)|\sqrt{b}-w^r(\frac{M_1(t)}{t})|+|w^r(\frac{M_1(t)}{t})-w(M_1(t),t))|\\
=&O(|\frac{\sqrt{b}t}{t}-\frac{M_1(t)}{t}|)+O(1+t)^{-\frac{3}{4}}\\
=&O(\frac{M_1(t)-\sqrt{b}t}{t})=O(\frac{(\ln(1+t)^{\beta})^{\frac{1}{2}}}{\sqrt{1+t}})
\end{aligned}
\end{equation}
When $\beta=\frac{1}{4}$, $t$ is suitably large, 

\begin{equation}\label{loc2}
\begin{aligned}
&\frac{a-v_-}{\sqrt{\pi}}\int_{(\ln(1+t)^{\beta})^{\frac{1}{2}}}^{+\infty}e^{-\xi^2}d\xi=\frac{a-v_-}{\sqrt{\pi}}\int_{(\ln(1+t)^{\frac{1}{4}})^{\frac{1}{2}}}^{+\infty}e^{-\xi^2}d\xi\\
> &\frac{a-v_-}{\sqrt{\pi}}\int_{(\ln(1+t)^{\frac{1}{4}})^{\frac{1}{2}}}^{(\ln(1+t)^{\frac{1}{3}})^{\frac{1}{2}}}e^{-\xi^2}d\xi\geq O(1+t)^{-\frac{1}{3}}>O(1+t)^{-\frac{1}{2}}\\
\end{aligned}
\end{equation}
Combining \eqref{loc1}-\eqref{loc2}, there exists a small constant $\beta>0$ such that \eqref{h1z1} holds when $t$ is suitably large.
 We complete the proof of \eqref{lem3.3}.

\end{proof}

\section{Stability analysis}

In order to prove theorem \ref{mt}, it suffice to show that the composite wave $(\hat{v},\hat{u})(x,t)$ constructed in Case 1 is stable. We consider the solutions $(v,u)(x,t)$ of \eqref{eq} in a neighborhood of $(\hat{v},\hat{u})(x,t)$, order the perturbation function $(\phi,\psi)(x,t)$ as 
\begin{equation}
(\phi,\psi)(x,t)=(v-\hat{v},u-\hat{u})(x,t).
\end{equation}
From \eqref{eq}-\eqref{inta} and \eqref{hatueq}-\eqref{hatin}, the Cauchy problem for $(\phi,\psi)(x,t)$ is equivalent to the following initial value problem:

\begin{equation}\label{phi-psi}
\left\{
\begin{aligned}
&\phi_t-\psi_x=-Q_1,\\
&\psi_t-(\sigma(v)-\sigma(\hat{v}))_x=\mu \psi_{xx}-H,\\
&(\phi,\psi)(x,0)=(\phi_0,\psi_0)(x)\rightarrow (0,0),\text{as}\quad x\rightarrow\pm\infty.
\end{aligned}
\right.
\end{equation}

Then we only need to show that $(\phi,\psi)(x,t)\rightarrow 0$ as $t\rightarrow\infty.$ Since the local existence theorem of the solution is standard, we only need an a priori estimate on the solution of \eqref{phi-psi}. Define the solution space for \eqref{phi-psi} by

\begin{equation}
\mathbb{X}(0,T)=\left\{
\begin{aligned}
(\phi,\psi)(x,t)\mid (\phi,\psi)(x,t)\in C(0,T;H^1); \phi_{x}(x,t)\in L^2(0,T;L^2); \psi_{x}(x,t)\in L^2(0,T;H^1).
\end{aligned}
\right\}.
\end{equation}
Order 
\begin{equation}\label{N(T)}
N(T)=\sup_{\tau\in[0,T]}\|(\phi,\psi)(\cdot,\tau)\|_{H^1(\mathbb{R})}\leq \epsilon_0,
\end{equation}

\begin{Lem}\label{4.1}Under assumptions of \eqref{t1} in Theorem \ref{mt}, if  $N(T)\leq\epsilon_0,$ then the following basic energy estimate holds.

\begin{equation}\label{ene1}
\begin{aligned}
&\sup_{t\in[0,T]}\|(\phi,\psi)(t)\|^2+G(t)+\int_0^t\|\psi_x\|^2dxd\tau\\
\leq &C(\|(\phi_0,\psi_0)\|^2+1)+\frac{1}{300}\int_0^t\|\phi_x\|^2d\tau+C\int_0^t\int_{-\infty}^{+\infty}\Xi_{2x}^2\phi^2dxd\tau
\end{aligned}
\end{equation}
where $G(t)$ is defined by

\begin{equation}\label{G(t)}
\begin{aligned}
G(t)=&-\int_0^t\int_{\hat{v}+\phi>a,\hat{v}>a}\hat{u}_{1x}\phi^2 dxd\tau-\int_0^t\int_{\hat{v}+\phi>a,-a<\hat{v}\leq a}\hat{u}_{1x}(\phi+\hat{v}-a)^2dxd\tau\\
&-\int_0^t\int_{-a<\hat{v}+\phi\leq a,\hat{v}>a}\hat{u}_{1x}(a-\hat{v})^2dxd\tau,
\end{aligned}
\end{equation}
and $-\hat{u}_{1x}=-u^r_x+\sqrt{b}\Xi_{2x}>0$.
\end{Lem}

\begin{proof}Multiplying $\eqref{phi-psi}_1$ by $(\sigma(v)-\sigma(\hat{v}))$ and $\eqref{phi-psi}_2$ by $\psi$ gives that

\begin{equation}\label{bs}
\begin{aligned}
&[\int_{\hat{v}}^{\hat{v}+\phi}\sigma(s)ds-\sigma(\hat{v})\phi]_t+[\frac{1}{2}\psi^2]_t-\hat{v}_t\bigg\{\sigma(\hat{v}+\phi)-\sigma(\hat{v})-\sigma'(\hat{v})\phi\bigg\}+(\cdots)_x\\
&+\mu\psi_x^2=(\sigma(\hat{v})-\sigma(v))Q_1-H\psi.
\end{aligned}
\end{equation}
Note that $\sigma'(v)>0$ for $v\in(-\infty,+\infty),$ this implies that 

\begin{equation}\label{taylor}
\int_{\hat{v}}^{\hat{v}+\phi}\sigma(s)ds-\sigma(\hat{v})\phi=F(v)-F(\hat{v})-F'(\hat{v})\phi=\frac{1}{2}F''(\theta\hat{v}+(1-\theta)v)\phi^2=\frac{1}{2}\sigma'(v_{\theta})\phi^2.
\end{equation}

Integrating \eqref{bs} over $[0,t]\times\mathbb{R}$ yields 

\begin{equation}\label{bse1}
\begin{aligned}
&\|(\phi,\psi)(t)\|^2+\int_0^t\int_{-\infty}^{+\infty}|\psi_x|^2dxd\tau\\
=&\|(\phi_0,\psi_0)\|^2+\int_0^t\int_{-\infty}^{+\infty}|Q_1\phi|-H\psi dxd\tau\\
&+\int_0^t\int_{-\infty}^{+\infty}\hat{v}_{\tau}\bigg\{\sigma(\hat{v}+\phi)-\sigma(\hat{v})-\sigma'(\hat{v})\phi\bigg\}dxd\tau.
\end{aligned}
\end{equation}
Estimates on terms in \eqref{bse1} will be separated as follows.

Similar as  \cite{HLM}, using \eqref{hatu}, \eqref{hatueq}, \eqref{HL1}, we have

\begin{equation}\label{bse2}
\begin{aligned}
\int_0^t\int_{-\infty}^{+\infty}H\psi dxd\tau
\leq &\int_0^t\|H\|_{L^1}\|\psi\|_{L^{\infty}}d\tau\\
\leq &\frac{1}{10}\int_0^t\int_{-\infty}^{+\infty}|\psi_x|^2dxd\tau+\frac{1}{10}\sup_{t\in[0,T]}\|\psi\|^2+C\int_0^t\|H\|_{L^1}^{\frac{4}{3}}d\tau
\end{aligned}
\end{equation}

\begin{equation}\label{bse3}
\begin{aligned}
&\int_0^t\int_{-\infty}^{+\infty}\hat{v}_{\tau}\bigg\{\sigma(\hat{v}+\phi)-\sigma(\hat{v})-\sigma'(\hat{v})\phi\bigg\}dxd\tau\\
\leq &\int_0^t\int_{-\infty}^{+\infty} \bigg[(u^r_x-\sqrt{b}\Xi_{2x})+(\Xi_{2xx},|Q_1|)\bigg]\bigg\{\sigma(\hat{v}+\phi)-\sigma(\hat{v})-\sigma'(\hat{v})\phi\bigg\}dxd\tau.
\end{aligned}
\end{equation}
Since $\hat{u}_{1x}:=(u^r-\sqrt{b}\Xi_{2})_x<0$, we have

\begin{equation}\label{a1}
\begin{aligned}
&\int_0^t\int_{-\infty}^{+\infty}\hat{u}_{1x}\bigg\{\sigma(\hat{v}+\phi)-\sigma(\hat{v})-\sigma'(\hat{v})\phi\bigg\} dxd\tau<0.
\end{aligned}
\end{equation}
In fact, the above integral \eqref{a1} can be
divided into four parts.

\begin{equation}\label{a21}
\begin{aligned}
&\int_0^t\int_{\hat{v}+\phi>a,\hat{v}>a}\hat{u}_{1x}\bigg\{\sigma(\hat{v}+\phi)-\sigma(\hat{v})-\sigma'(\hat{v})\phi\bigg\}dxd\tau\\
\leq &C^{-1}\int_0^t\int_{\hat{v}+\phi>a,\hat{v}>a}\hat{u}_{1x}\phi^2 dxd\tau,
\end{aligned}
\end{equation}

\begin{equation}\label{a22}
\begin{aligned}
&\int_0^t\int_{\hat{v}+\phi>a,-a<\hat{v}\leq a}\hat{u}_{1x}\bigg\{\sigma(\hat{v}+\phi)-\sigma(\hat{v})-\sigma'(\hat{v})\phi\bigg\}dxd\tau\\
\leq &C^{-1}\int_0^t\int_{\hat{v}+\phi>a,-a<\hat{v}\leq a}\hat{u}_{1x}(\phi+\hat{v}-a)^2dxd\tau,
\end{aligned}
\end{equation}

\begin{equation}\label{a23}
\begin{aligned}
&\int_0^t\int_{-a<\hat{v}+\phi\leq a,\hat{v}>a}\hat{u}_{1x}\bigg\{\sigma(\hat{v}+\phi)-\sigma(\hat{v})-\sigma'(\hat{v})\phi\bigg\}dxd\tau\\
\leq & C^{-1}\int_0^t\int_{-a<\hat{v}+\phi\leq a,\hat{v}>a}\hat{u}_{1x}(a-\hat{v})^2dxd\tau,
\end{aligned}
\end{equation}

\begin{equation}\label{a24}
\begin{aligned}
\int_0^t\int_{-a<\hat{v}+\phi\leq a,-a<\hat{v}\leq a}\hat{u}_{1x}\bigg\{\sigma(\hat{v}+\phi)-\sigma(\hat{v})-\sigma'(\hat{v})\phi\bigg\}dxd\tau=0.
\end{aligned}
\end{equation}
As for other terms in \eqref{bse3}, we estimate them as follows:

\begin{equation}\label{hk}
\begin{aligned}
&\int_0^t\int_{-\infty}^{+\infty}\Xi_{2xx}\bigg\{\sigma(\hat{v}+\phi)-\sigma(\hat{v})-\sigma'(\hat{v})\phi\bigg\}dxd\tau\\
=&\int_0^t\int_{-\infty}^{Z_1(\tau)}\Xi_{2xx}\bigg\{\sigma(\hat{v}+\phi)-\sigma(\hat{v})-\sigma'(\hat{v})\phi\bigg\}dxd\tau\\
&+\int_0^t\int_{Z_1(\tau)}^{X_1(\tau)}\Xi_{2xx}\bigg\{\sigma(\hat{v}+\phi)-\sigma(\hat{v})-\sigma'(\hat{v})\phi\bigg\}dxd\tau\\
&+\int_0^t\int_{X_1(\tau)}^{+\infty}\Xi_{2xx}\bigg\{\sigma(\hat{v}+\phi)-\sigma(\hat{v})-\sigma'(\hat{v})\phi\bigg\}dxd\tau.
\end{aligned}
\end{equation}
Since when $x\geq X_1(t),\Xi_{2xx}<0$(see \eqref{Xi2x}), by using the analyse in  \eqref{a21}-\eqref{a24},

\begin{equation*}
0\leq \sigma(\hat{v}+\phi)-\sigma({\hat{v}})-\sigma'(\hat{v})\phi\leq C\phi^2,
\end{equation*}
then 
\begin{equation}\label{hk1}
\begin{aligned}
\int_0^t\int_{X_1(\tau)}^{+\infty}\Xi_{2xx}\bigg\{\sigma(\hat{v}+\phi)-\sigma({\hat{v}})-\sigma'(\hat{v})\phi\bigg\} dxd\tau<0.
\end{aligned}
\end{equation}
We omit the estimation of the third term in \eqref{hk}. Similar as \eqref{Z1-3},

\begin{equation}\label{hk2}
\begin{aligned}
&\int_0^t\int_{Z_1(\tau)}^{X_1(\tau)}\Xi_{2xx}\bigg\{\sigma(\hat{v}+\phi)-\sigma({\hat{v}})-\sigma'(\hat{v})\phi\bigg\} dxd\tau\\
\leq &C\int_0^t\|\phi\|_{L^{\infty}}^2d\tau \int_{Z_1(\tau)}^{X_1(\tau)}|\Xi_{2xx}|dx\\
\leq &C\delta\int_0^t(1+\tau)^{-1}\|\phi\|\|\phi_x\|d\tau\\
\leq &C\delta\int_0^t\|\phi_x(\tau)\|^2d\tau+C\delta\int_0^t(1+\tau)^{-2}\|\phi(\tau)\|^2d\tau.
\end{aligned}
\end{equation}
As for $x\leq Z_1(t),$ we have obtained that $\hat{v}\leq a$ by \eqref{Z1}, then $\sigma'(\hat{v})=\sigma'(a)$, integration by parts, 
\begin{equation}\label{herk}
\begin{aligned}
&\int_0^t\int_{-\infty}^{Z_1(\tau)}\Xi_{2xx}\bigg\{\sigma(\hat{v}+\phi)-\sigma({\hat{v}})-\sigma'(\hat{v})\phi\bigg\}dxd\tau\\
=&\int_0^t\int_{-\infty}^{Z_1(\tau)}\Xi_{2xx}\bigg\{\sigma(\hat{v}+\phi)-\sigma({\hat{v}})-\sigma'(a)\phi\bigg\}dxd\tau\\
\leq & C\int_0^t\big(\Xi_{2x}\phi^2\big)(Z_1(\tau),\tau)d\tau-\int_0^t\int_{-\infty}^{Z_1(\tau)}\Xi_{2x}\bigg\{\sigma(\hat{v}+\phi)-\sigma({\hat{v}})-\sigma'(a)\phi\bigg\}_x dxd\tau
\end{aligned}
\end{equation}
Here by using \eqref{lem3.3},

\begin{equation}\label{herk1}
\begin{aligned}
&\int_0^t\big(\Xi_{2x}\phi^2\big)(Z_1(\tau),\tau)d\tau\\
\leq &C\delta \int_0^t\|\phi\|^2(1+\tau)^{-1}e^{-\frac{(Z_1(\tau)-\sqrt{b}\tau)^2}{2\mu(1+\tau)}}d\tau+C\delta\int_0^t\|\phi_x(\tau)\|^2d\tau\\
\leq &C\delta\int_0^t\|\phi\|^2(1+\tau)^{-1-\beta}d\tau+C\delta\int_0^t\|\phi_x(\tau)\|^2d\tau\\
\leq & C\delta \bigg(\|\phi\|^2+\int_0^t\|\phi_x(\tau)\|^2d\tau\bigg).
\end{aligned}
\end{equation}
As for another term in \eqref{herk}, 
\begin{equation}\label{herk2}
\begin{aligned}
&-\int_0^t\int_{-\infty}^{Z_1(\tau)}\Xi_{2x}\bigg\{\sigma(\hat{v}+\phi)-\sigma({\hat{v}})-\sigma'(a)\phi\bigg\}_x dxd\tau\\
=&-\int_0^t\int_{-\infty}^{Z_1(\tau)}\Xi_{2x}\bigg\{\big[\sigma'(\hat{v}+\phi)-\sigma'(\hat{v})\big](\phi_x+\hat{v}_x)\bigg\}dxd\tau\\
\leq &\frac{1}{400}\int_0^t\|\phi_x(\tau)\|^2d\tau+C\int_0^t\int_{-\infty}^{Z_1(\tau)}\Xi_{2x}^2\phi^2dxd\tau-\int_0^t\int_{-\infty}^{Z_1(\tau)}\Xi_{2x}\big[\sigma'(\hat{v}+\phi)-\sigma'(\hat{v})\big]\hat{v}_x dxd\tau,
\end{aligned}
\end{equation}

\begin{equation}\label{herk3}
\begin{aligned}
&-\int_0^t\int_{-\infty}^{Z_1(\tau)}\Xi_{2x}\big[\sigma'(\hat{v}+\phi)-\sigma'(\hat{v})\big]\hat{v}_x dxd\tau\\
=&-C\int_0^t\int_{\hat{v}+\phi>a,-a<\hat{v}\leq a}\Xi_{2x}[(\phi+\hat{v}-a)](v^r_x+\Xi_{2x}+\Xi_{2xx})dxd\tau
\end{aligned}
\end{equation}

\begin{equation}\label{herk4}
-\int_0^t\int_{\hat{v}+\phi>a,-a<\hat{v}\leq a}\Xi_{2x}[(\phi+\hat{v}-a)](v^r_x+\Xi_{2x})dxd\tau<0
\end{equation}

\begin{equation}\label{herk5}
\begin{aligned}
&\int_0^t\int_{\hat{v}+\phi>a,-a<\hat{v}\leq a}\Xi_{2x}[(\phi+\hat{v}-a)]\Xi_{2xx}dxd\tau \\
\leq &C\int_0^t\int_{-\infty}^{+\infty}\Xi_{2x}^2\phi^2dxd\tau+\int_0^t\int_{-\infty}^{+\infty}\Xi_{2xx}^2dxd\tau\\
\leq &C\int_0^t\int_{-\infty}^{+\infty}\Xi_{2x}^2\phi^2dxd\tau+C\delta\int_0^t\int_{-\infty}^{+\infty}(1+\tau)^{-2}e^{-\frac{c(x-\sqrt{b}\tau)^2}{1+\tau}} dxd\tau
\end{aligned}
\end{equation}
Inserting \eqref{herk1}-\eqref{herk5} into \eqref{herk}, we could obtain

\begin{equation}\label{hk3}
\begin{aligned}
&\int_0^t\int_{-\infty}^{Z_1(\tau)}\Xi_{2xx}\bigg\{\sigma(\hat{v}+\phi)-\sigma({\hat{v}})-\sigma'(\hat{v})\phi\bigg\}dxd\tau\\
\leq &C\delta (1+\|\phi\|^2)+\frac{1}{300}\int_0^t\int_0^t\|\phi_x(\tau)\|^2d\tau+C\int_0^t\int_{-\infty}^{+\infty}\Xi_{2x}^2\phi^2dxd\tau.
\end{aligned}
\end{equation}


\begin{equation}\label{Qphi}
\begin{aligned}
\int_0^t\int_{-\infty}^{+\infty}|Q_1\phi|dxd\tau\leq O(\delta)\int_0^t(1+\tau)^{-1}\|\phi\|^{\frac{1}{2}}\|\phi_x\|^{\frac{1}{2}}d\tau\leq O(\delta)\bigg(1+\sup_{t\in[0,T]}\|\phi\|^2+\int_0^t\|\phi_x\|^2d\tau\bigg).
\end{aligned}
\end{equation}
Combining  \eqref{bse2}-\eqref{Qphi}, we complete the proof of \eqref{ene1}.
\end{proof}

\begin{Lem}(heat kernal estimates)\label{4.2}Under assumptions of \eqref{t1} in Theorem \ref{mt}, if  $N(T)\leq\epsilon_0,$ then the following energy estimate holds.

\begin{equation}\label{hke}
\begin{aligned}
&\int_0^t\int_{\mathbb{R}}w^2(\phi^2+\psi^2)dxd\tau\\
\leq &C\int_0^t\|(\phi_{x},\psi_{x})\|^2d\tau+C(\|(\phi,\psi)\|^2+1)+CG(t)\\
&+C(\delta+N(T))\int_0^t\int_{\mathbb{R}}w\bigg(\sigma(v)-\sigma(\hat{v})-\sigma'(\hat{v})\phi \bigg) dx d\tau.
\end{aligned}
\end{equation}
where heat-kernal function $w(x,t)$ is defined as

\begin{equation}\label{wx}
\begin{aligned}
w(x,t)=\frac{1}{\sqrt{2\mu(1+t)}}e^{-\frac{(x-\sqrt{b}t)^2}{2\mu(1+t)}}
\end{aligned}
\end{equation}

\end{Lem}

\begin{proof}Order $\eta=x-\sqrt{b}t,$ then \eqref{hke} becomes

\begin{equation}\label{hke1}
\begin{aligned}
&\int_0^t\int_{\mathbb{R}}w^2(\eta,\tau)(\phi^2+\psi^2)d\eta d\tau\\
\leq &C\int_0^t\|(\phi_{\eta},\psi_{\eta})\|^2d\tau+C(\|(\phi,\psi)\|^2+1)+CG(t)\\
&+C(\delta+N(T))\int_0^t\int_{\mathbb{R}}w\bigg(\sigma(v)-\sigma(\hat{v})-\sigma'(\hat{v})\phi \bigg) d\eta d\tau.
\end{aligned}
\end{equation}
where $w(\eta,t)$ is defined as 

\begin{equation}\label{weta}
\begin{aligned}
w(\eta,t)=\frac{1}{\sqrt{2\mu(1+t)}}e^{-\frac{\eta^2}{2\mu(1+t)}}.
\end{aligned}
\end{equation}
To prove \eqref{hke1}, we divided \eqref{hke1} into two parts;

\begin{equation}\label{hke1-1}
\begin{aligned}
&\int_0^t\int_{\mathbb{R}}w^2(\eta,\tau)(\sqrt{b}\phi-\psi)^2d\eta d\tau\\
\leq &C\|(\phi,\psi)\|^2+C\int_0^t\|(\phi_\eta,\psi_\eta)\|^2 d\tau+\int_0^t\bigg<(\sqrt{b}\phi-\psi)_{\tau},(\sqrt{b}\phi-\psi)g^2\bigg>d\tau\\
\leq &C\int_0^t\|(\phi_{\eta},\psi_{\eta})\|^2d\tau+\frac{1}{80}\int_0^t\int_{-\infty}^{+\infty}w^2(\eta,\tau)(\phi^2+\psi^2)d\eta d\tau+C(\|(\phi,\psi)\|^2+1)\\
&+CG(t)+C(\delta+N(T))\int_0^t\int_{\mathbb{R}}w\bigg(\sigma(v)-\sigma(\hat{v})-\sigma'(\hat{v})\phi \bigg) d\eta d\tau.
\end{aligned}
\end{equation}
where we use the method introduced in \cite{HLM}. Another inequality is

\begin{equation}\label{hke1-2}
\begin{aligned}
&\int_0^t\int_{\mathbb{R}}w^2(\eta,\tau)(\sqrt{b}\phi+\psi)^2d\eta d\tau\\
\leq &C\int_0^t\|(\phi_{\eta},\psi_{\eta})\|^2d\tau+\frac{1}{80}\int_0^t\int_{-\infty}^{+\infty}w^2(\eta,\tau)(\phi^2+\psi^2)d\eta d\tau+C(\|(\phi,\psi)\|^2+1)\\
&+CG(t)+C(\delta+N(T))\int_0^t\int_{\mathbb{R}}w\bigg(\sigma(v)-\sigma(\hat{v})-\sigma'(\hat{v})\phi \bigg) d\eta d\tau.
\end{aligned}
\end{equation}

\noindent{Step 1:} Doing the coordinates transformation $(x,t)\rightarrow(\eta,t),$ then \eqref{phi-psi} turns to

\begin{equation}\label{pbs}
\left\{
\begin{aligned}
&\phi_t-\sqrt{b}\phi_{\eta}-\psi_\eta=-Q_1,\\
&\psi_t-\sqrt{b}\psi_{\eta}-(\sigma(v)-\sigma(\hat{v}))_{\eta}=\mu \psi_{\eta\eta}-H,\\
&(\phi,\psi)(\eta,0)=(\phi_0,\psi_0)(\eta)\rightarrow (0,0),\text{as}\quad \eta\rightarrow\pm\infty.
\end{aligned}
\right.
\end{equation}

$\eqref{pbs}_1\times \sqrt{b}-\eqref{pbs}_2$ gives that

\begin{equation}\label{bphi-si}
\begin{aligned}
(\sqrt{b}\phi-\psi)_t+\bigg(\sigma(v)-\sigma(\hat{v})-b\phi\bigg)_{\eta}=-\sqrt{b}Q_1-\mu\psi_{\eta\eta}+H.
\end{aligned}
\end{equation}
Order $h=\sqrt{b}\phi-\psi$ and remember $w(\eta,t)$ in \eqref{weta}, we define new functions
\begin{equation}\label{f}
\begin{aligned}
&f(\eta,t)=\int_{-\infty}^{\eta}w^2(y,t)dy,\quad f_{\eta}(\eta,t)=w^2(\eta,t),
\end{aligned}
\end{equation}

\begin{equation}\label{g}
\begin{aligned}
&g=\int_{-\infty}^{\eta}w(y,\tau)dy,\quad \|g\|_{L^{\infty}}\leq C,\\
&\|g_{\eta}\|_{L^{\infty}}\leq C(1+t)^{-\frac{1}{2}},\quad \|g_t\|_{L^{\infty}}\leq  C(1+t)^{-1}.
\end{aligned}
\end{equation}
Using Huang-Matsumura inequality, \cite{HLM},
we get

\begin{equation}\label{bh-s1}
\begin{aligned}
&\int_0^t\int_{\mathbb{R}}(\sqrt{b}\phi-\psi)^2w^2(\eta,\tau)d\eta d\tau\\
\leq &C\|(\phi_0,\psi_0)\|^2+C\int_0^t\|(\phi_{\eta},\psi_{\eta})\|^2d\tau+C\int_0^t\bigg<(\sqrt{b}\phi-\psi)_{\tau},(\sqrt{b}\phi-\psi)g^2\bigg>d\tau.
\end{aligned}
\end{equation}
Here 

\begin{equation}\label{bh-s2}
\begin{aligned}
&\int_0^t\bigg<(\sqrt{b}\phi-\psi)_{\tau},(\sqrt{b}\phi-\psi)g^2\bigg>d\tau\\
=&-\int_0^t\int_{\mathbb{R}}\bigg(\sigma(v)-\sigma(\hat{v})-b\phi\bigg)_{\eta}(\sqrt{b}\phi-\psi)g^2 d\eta d\tau\\
&+\int_0^t\int_{\mathbb{R}}\bigg[-\sqrt{b}Q_1-\mu\psi_{\eta\eta}+H\bigg](\sqrt{b}\phi-\psi)g^2 d\eta d\tau.
\end{aligned}
\end{equation}
We estimate \eqref{bh-s2} item by item.

\begin{equation}\label{J}
\begin{aligned}
&-\int_0^t\int_{\mathbb{R}}\bigg(\sigma(v)-\sigma(\hat{v})-b\phi\bigg)_{\eta}(\sqrt{b}\phi-\psi)g^2 d\eta d\tau\\
=&\int_0^t\int_{\mathbb{R}}\bigg(\sigma(v)-\sigma(\hat{v})-b\phi \bigg)(\sqrt{b}\phi_{\eta}-\psi_{\eta})g^2 d\eta d\tau\\
&+\int_0^t\int_{\mathbb{R}}\bigg(\sigma(v)-\sigma(\hat{v})-b\phi \bigg)(\sqrt{b}\phi-\psi)2gg_{\eta}d\eta d\tau\\
=&J_1+J_2.
\end{aligned}
\end{equation}

\begin{equation}\label{J1}
\begin{aligned}
J_1=&\int_0^t\int_{\mathbb{R}}\bigg(\sigma(v)-\sigma(\hat{v})-b\phi \bigg)(\sqrt{b}\phi_{\eta}-\psi_{\eta})g^2 d\eta d\tau\\
=&\int_0^t\int_{\mathbb{R}}\bigg(\sigma(v)-\sigma(\hat{v})-b\phi \bigg)(2\sqrt{b}\phi_{\eta}-\phi_{\tau}-Q_1)g^2 d\eta d\tau\\
=&\int_0^t\int_{\mathbb{R}}2\sqrt{b}\bigg[\int_{\hat{v}}^{\hat{v}+\phi}\sigma(s)ds-\sigma(\hat{v})\phi-\frac{1}{2}b\phi^2\bigg]_{\eta}g^2 d\eta d\tau\\
&-\int_0^t\int_{\mathbb{R}}\bigg[\int_{\hat{v}}^{\hat{v}+\phi}\sigma(s)ds-\sigma(\hat{v})\phi-\frac{1}{2}b\phi^2\bigg]_{\tau}g^2d\eta d\tau\\
&-\int_0^t\int_{\mathbb{R}}(2\sqrt{b}\hat{v}_{\eta}-\hat{v}_{\tau})\bigg(\sigma(v)-\sigma(\hat{v})-\sigma'(\hat{v})\phi \bigg)g^2 d\eta d\tau\\
&-\int_0^t\int_{\mathbb{R}}\bigg(\sigma(v)-\sigma(\hat{v})-b\phi \bigg)Q_1g^2 d\eta d\tau\\
=&J_{11}+J_{12}+J_{13}+J_{14}.
\end{aligned}
\end{equation}
Similar as \eqref{bse3}, we have
\begin{equation}
\begin{aligned}
J_{13}+J_{14}\leq &CG(t)+O(\delta)\bigg(1+\sup_{t\in[0,T]}\|\phi\|^2+\int_0^t\|\phi_{\eta}\|^2d\tau\bigg)\\
&+C\int_0^t\int_{-\infty}^{+\infty}\Xi_{2\eta}^2\phi^2d\eta d\tau.
\end{aligned}
\end{equation}
Compared with \eqref{weta}, the relationship between $\Xi_{2\eta}$ and $w(\eta,t)$ is 

\begin{equation}\label{Xiw}
\begin{aligned}
\Xi_{2\eta}(\eta,t)=\frac{(a-v_-)}{\sqrt{2\mu(1+t)}}e^{-\frac{\eta^2}{2\mu(1+t)}}=O(\delta)w(\eta,t).
\end{aligned}
\end{equation}
Then it yields that

\begin{equation}\label{J134}
\begin{aligned}
J_{13}+J_{14}\leq &CG(t)+O(\delta)\bigg(1+\sup_{t\in[0,T]}\|\phi\|^2+\int_0^t\|\phi_{\eta}\|^2d\tau\bigg)\\
&+O(\delta)\int_0^t\int_{-\infty}^{+\infty}w^2\phi^2d\eta d\tau.
\end{aligned}
\end{equation}
As for $J_{11},$ after integration by parts with $\eta,$ it yields that
\begin{equation}\label{J11}
\begin{aligned}
J_{11}=&4\sqrt{b}\int_0^t\int_{\mathbb{R}}\bigg[\int_{\hat{v}}^{\hat{v}+\phi}\sigma(s)ds-\sigma(\hat{v})\phi-\frac{1}{2}b\phi^2\bigg]gg_{\eta} d\eta d\tau\\
=&4\sqrt{b}\int_0^t\int_{\hat{v}+\phi>a,\hat{v}>a}\bigg[\int_{\hat{v}}^{\hat{v}+\phi}\sigma(s)ds-\sigma(\hat{v})\phi-\frac{1}{2}b\phi^2\bigg]gg_{\eta}d\eta d\tau\\
&+4\sqrt{b}\int_0^t\int_{\hat{v}+\phi>a,-a<\hat{v}\leq a}\bigg[\int_{\hat{v}}^{\hat{v}+\phi}\sigma(s)ds-\sigma(\hat{v})\phi-\frac{1}{2}b\phi^2\bigg]gg_{\eta}d\eta d\tau\\
&+4\sqrt{b}\int_0^t\int_{-a<\hat{v}+\phi\leq a,\hat{v}>a}\bigg[\int_{\hat{v}}^{\hat{v}+\phi}\sigma(s)ds-\sigma(\hat{v})\phi-\frac{1}{2}b\phi^2\bigg]gg_{\eta}d\eta d\tau\\
&+4\sqrt{b}\int_0^t\int_{-a<\hat{v}+\phi\leq a,-a<\hat{v}\leq a}\bigg[\int_{\hat{v}}^{\hat{v}+\phi}\sigma(s)ds-\sigma(\hat{v})\phi-\frac{1}{2}b\phi^2\bigg]gg_{\eta}d\eta d\tau\\
=& J_{111}+J_{112}+J_{113}+J_{114}.
\end{aligned}
\end{equation}
Using the Taylor expansion in \eqref{taylor}, when $\hat{v}+\phi>a,\hat{v}>a,$
\begin{equation}
\begin{aligned}
&\int_{\hat{v}}^{\hat{v}+\phi}\sigma(s)ds-\sigma(\hat{v})\phi-\frac{1}{2}b\phi^2,\\
=&\frac{1}{2}\sigma'(\hat{v})\phi^2+\frac{1}{3!}\sigma''(\hat{v})\phi^3+O(1)\phi^4-\frac{1}{2}b\phi^2,
\end{aligned}
\end{equation}

\begin{equation}\label{J111}
\begin{aligned}
J_{111}=&\int_0^t\int_{\hat{v}+\phi>a,\hat{v}>a}\bigg[\int_{\hat{v}}^{\hat{v}+\phi}\sigma(s)ds-\sigma(\hat{v})\phi-\frac{1}{2}b\phi^2\bigg]gg_{\eta} d\eta d\tau\\
\leq &C\int_0^t\int_{\hat{v}+\phi>a,\hat{v}>a}\bigg\{\frac{1}{2}(\sigma'(\hat{v})-\sigma'(a))\phi^2+\frac{1}{6}\sigma''(\hat{v})\phi^3+\phi^4\bigg\}w(\eta,\tau) d\eta d\tau\\
\leq  &C\|(\hat{v}-a,\phi)\|_{\infty}\int_0^t\int_{\hat{v}+\phi>a,\hat{v}>a}w\phi^2d\eta d\tau,
\end{aligned}
\end{equation}
when $\hat{v}+\phi>a,-a<\hat{v}\leq a,$ from \eqref{sigma-1}\ $\sigma(\hat{v})=b\hat{v},$ then

\begin{equation}
\begin{aligned}
&\int_{\hat{v}}^{\hat{v}+\phi}\sigma(s)ds-\sigma(\hat{v})\phi-\frac{1}{2}b\phi^2\\
=&\int_{\hat{v}}^{a}bs ds+\int_{a}^{\hat{v}+\phi}\sigma(s)ds-b\hat{v}\phi-\frac{1}{2}b\phi^2\\
=&\frac{1}{2}ba^2-\frac{1}{2}b\hat{v}^2+ba(\hat{v}+\phi-a)+\frac{1}{2}\sigma'(a)(\hat{v}+\phi-a)^2+O(\hat{v}+\phi-a)^3-b\hat{v}\phi-\frac{1}{2}b\phi^2\\
=&O(\hat{v}+\phi-a)^3,
\end{aligned}
\end{equation}

\begin{equation}\label{J112}
\begin{aligned}
J_{112}=&\int_0^t\int_{\hat{v}+\phi>a,-a<\hat{v}\leq a}\bigg[\int_{\hat{v}}^{\hat{v}+\phi}\sigma(s)ds-\sigma(\hat{v})\phi-\frac{1}{2}b\phi^2\bigg]gg_{\eta} d\eta d\tau\\
\leq &C\int_0^t\int_{\hat{v}+\phi>a,-a<\hat{v}\leq a}w(\hat{v}+\phi-a)^3d\eta d\tau\\
\leq &C\|\hat{v}+\phi-a\|_{\infty}\int_0^t\int_{\hat{v}+\phi>a,-a<\hat{v}\leq a}w(\hat{v}+\phi-a)^2d\eta d\tau
\end{aligned}
\end{equation}
when $-a<\hat{v}+\phi\leq a,\hat{v}>a,$

\begin{equation}
\begin{aligned}
&\int_{\hat{v}}^{\hat{v}+\phi}\sigma(s)ds-\sigma(\hat{v})\phi-\frac{1}{2}b\phi^2\\
=&\int_{\hat{v}}^{a}\sigma(s) ds+\int_{a}^{\hat{v}+\phi}bsds-\sigma(\hat{v})\phi-\frac{1}{2}b\phi^2\\
=&\frac{1}{2}b\hat{v}^2+b\hat{v}\phi-\frac{1}{2}ba^2-\bigg\{ba(\hat{v}-a)+\frac{1}{2}b(\hat{v}-a)^2+O(\hat{v}-a)^3+\sigma(\hat{v})\phi\bigg\}\\
=&O(\hat{v}-a)^3-\bigg(\sigma(\hat{v})-\sigma(a)-\sigma'(a)(\hat{v}-a)\bigg)\phi
,\end{aligned}
\end{equation}

\begin{equation}\label{J113}
\begin{aligned}
J_{113}=&\int_0^t\int_{-a<\hat{v}+\phi\leq a,\hat{v}>a}\bigg[\int_{\hat{v}}^{\hat{v}+\phi}\sigma(s)ds-\sigma(\hat{v})\phi-\frac{1}{2}b\phi^2\bigg]gg_{\eta}d\eta d\tau \\
\leq &\int_0^t\int_{-a<\hat{v}+\phi\leq a,\hat{v}>a}w\bigg[(\hat{v}-a)^3+(\hat{v}-a)^2\phi\bigg]d\eta d\tau\\
\leq &\|(\hat{v}-a,\phi)\|_{\infty}\int_0^t\int_{-a<\hat{v}+\phi\leq a,\hat{v}>a}w(\hat{v}-a)^2 d\eta d\tau
\end{aligned}
\end{equation}
when $-a<\hat{v}+\phi\leq a,-a<\hat{v}\leq a,$

\begin{equation}
\begin{aligned}
&\int_{\hat{v}}^{\hat{v}+\phi}\sigma(s)ds-\sigma(\hat{v})\phi-\frac{1}{2}b\phi^2\\
=&\int_{\hat{v}}^{\hat{v}+\phi}bsds-b\hat{v}\phi-\frac{1}{2}b\phi^2\\
=&\frac{1}{2}b(\hat{v}+\phi)^2-\frac{1}{2}b\hat{v}^2-b\hat{v}\phi-\frac{1}{2}b\phi^2\\=&0,
\end{aligned}
\end{equation}

\begin{equation}\label{J114}
\begin{aligned}
J_{114}=\int_0^t\int_{-a<\hat{v}+\phi\leq a,-a<\hat{v}\leq a}\bigg[\int_{\hat{v}}^{\hat{v}+\phi}\sigma(s)ds-\sigma(\hat{v})\phi-\frac{1}{2}b\phi^2\bigg]gg_{\eta}d\eta d\tau=0.
\end{aligned}
\end{equation}
Compared with \eqref{a21}-\eqref{a24}, the analysis in \eqref{J111},\eqref{J112},\eqref{J113},\eqref{J114} implies that

\begin{equation}\label{J11f}
\begin{aligned}
J_{11}\leq &C(\delta+N(T))\int_0^t\int_{\mathbb{R}}w(\eta,\tau)\bigg(\sigma(v)-\sigma(\hat{v})-\sigma'(\hat{v})\phi \bigg)d\eta d\tau.
\end{aligned}
\end{equation}
Also, 

\begin{equation}\label{J12}
\begin{aligned}
J_{12}=&\int_0^t\int_{\mathbb{R}}\bigg[\bigg(\int_{\hat{v}}^{\hat{v}+\phi}\sigma(s)ds-\sigma(\hat{v})\phi-\frac{1}{2}b\phi^2\bigg)g^2\bigg]_{\tau}d\eta d\tau\\
&-\int_0^t\int_{\mathbb{R}}\bigg[\int_{\hat{v}}^{\hat{v}+\phi}\sigma(s)ds-\sigma(\hat{v})\phi-\frac{1}{2}b\phi^2\bigg]2gg_{\tau} d\eta d\tau\\
\leq &C\|\phi\|^2+CJ_{11}.
\end{aligned}
\end{equation}
Then from \eqref{J134}, \eqref{J11f},\eqref{J12}, we complete the estimation of $J_1$ in \eqref{J1}.

Next, for $J_2$ in \eqref{J}, we see that

\begin{equation}\label{J2}
\begin{aligned}
J_2=&\int_0^t\int_{\mathbb{R}}\bigg(\sigma(v)-\sigma(\hat{v})-b\phi \bigg)(\sqrt{b}\phi-\psi)2gg_{\eta}d\eta d\tau\\
\leq &\int_0^t\int_{\mathbb{R}}w \bigg(\sigma(v)-\sigma(\hat{v})-\sigma'(\hat{v})\phi \bigg)(\sqrt{b}\phi-\psi)d\eta d\tau\\
&+\int_0^t\int_{\mathbb{R}}w\bigg(\sigma'(\hat{v})\phi-b\phi\bigg) (\sqrt{b}\phi-\psi)d\eta d\tau\\
\leq &\|(\phi,\psi)\|_{\infty}\int_0^t\int_{\mathbb{R}}w\bigg(\sigma(v)-\sigma(\hat{v})-\sigma'(\hat{v})\phi \bigg)d\eta d\tau\\
&+\int_0^t\int_{\mathbb{R}}w\bigg(\sigma'(\hat{v})-\sigma'(a)\bigg)\phi (\sqrt{b}\phi-\psi)d\eta d\tau=J_{21}+J_{22},
\end{aligned}
\end{equation}
and
\begin{equation}\label{J22}
\begin{aligned}
J_{22}=\int_0^t\int_{\mathbb{R}}w\bigg(\sigma'(\hat{v})-\sigma'(a)\bigg)\phi (\sqrt{b}\phi-\psi)d\eta d\tau
=\int_0^t\int_{\hat{v}>a}w\bigg(\sigma'(\hat{v})-\sigma'(a)\bigg)\phi (\sqrt{b}\phi-\psi)d\eta d\tau.
\end{aligned}
\end{equation}
When $\hat{v}>a\geq \hat{v}_1,$ obviously, $|\hat{v}-a|\leq |\hat{v}-\hat{v}_1|\leq \frac{\mu}{4\sqrt{b}}\Xi_{2\eta}.$  By using \eqref{Xiw},

\begin{equation}\label{J221}
\begin{aligned}
J_{221}=&\int_0^t\int_{\hat{v}>a\geq \hat{v}_1}w\bigg(\sigma'(\hat{v})-\sigma'(a)\bigg)\phi (\sqrt{b}\phi-\psi)d\eta d\tau\\
\leq &C\int_0^t\int_{\hat{v}>a\geq \hat{v}_1}w|\hat{v}-a|(\phi^2+\psi^2)d\eta d\tau\\
\leq &C\int_0^t\int_{\hat{v}>a\geq \hat{v}_1}w\Xi_{2\eta}(\phi^2+\psi^2)d\eta d\tau\\
\leq &C\delta\int_0^t\int_{\hat{v}>a\geq \hat{v}_1} w^2(\phi^2+\psi^2)d\eta d\tau.
\end{aligned}
\end{equation}
When $\hat{v}\geq \hat{v}_1>a,$ it implies that $x\geq X_1(t),$ then we have
\begin{equation}\label{J222}
\begin{aligned}
J_{222}=&\int_0^t\int_{\hat{v}\geq \hat{v}_1>a}w\bigg(\sigma'(\hat{v})-\sigma'(\hat{v}_1)+\sigma'(\hat{v}_1)-\sigma'(a)\bigg)\phi (\sqrt{b}\phi-\psi)d\eta d\tau\\
\leq &C\delta\int_0^t\int_{\hat{v}\geq \hat{v}_1>a} w^2(\phi^2+\psi^2)d\eta d\tau+C\int_0^t\int_{\hat{v}\geq \hat{v}_1>a}w(\hat{v}_1-a)(\phi^2-\phi\psi)d\eta d\tau.
\end{aligned}
\end{equation}
In detail, by \eqref{l3}, if $x\geq X_1(t)$, that is $\eta\geq X_1(t)-\sqrt{b}t>0$, we have

\begin{equation}
\begin{aligned}
&\int_0^t\int_{\hat{v}_1>a}w|\hat{v}_1-a|(\phi^2-\phi\psi)d\eta d\tau\\
\leq &C\|\hat{v}_1-a\|_{\infty}\int_0^t(1+\tau)^{-\frac{1}{2}}\int_{X_1(\tau)-\sqrt{b}\tau}^{+\infty}(\phi^2+\psi^2)e^{-\frac{\eta^2}{2\mu(1+\tau)}}d\eta d\tau\\
\leq &C\delta\int_0^t(1+\tau)^{-\frac{1}{2}}\int_{X_1(\tau)-\sqrt{b}\tau}^{Y_1(\tau)-\sqrt{b}\tau}(\phi^2+\psi^2)e^{-\frac{\eta^2}{2\mu(1+\tau)}}d\eta d\tau\\
&+C\delta\int_0^t(1+\tau)^{-\frac{1}{2}}\int_{Y_1(\tau)-\sqrt{b}\tau}^{+\infty}(\phi^2+\psi^2)e^{-\frac{\eta^2}{2\mu(1+\tau)}}d\eta d\tau,
\end{aligned}
\end{equation}
where 

\begin{equation}
\begin{aligned}
Y_1(t)=\sqrt{b}(1+t)+(C+\ln(1+t)^{\frac{1}{2(1+\epsilon)}+\epsilon})^{\frac{1}{2}}\sqrt{2\mu(1+t)}\geq X_1(t),
\end{aligned}
\end{equation}
compared with \eqref{l3}. When $\epsilon>0$ is suitably small, $\frac{1}{2(1+\epsilon)}+\epsilon>\frac{1}{2},$ it yields that 

\begin{equation}
\begin{aligned}
&C\delta\int_0^t(1+\tau)^{-\frac{1}{2}}\int_{Y_1(\tau)-\sqrt{b}\tau}^{+\infty}(\phi^2+\psi^2)e^{-\frac{\eta^2}{2\mu(1+\tau)}}d\eta d\tau\\
\leq &C\delta\int_0^t(1+\tau)^{-\frac{1}{2}-\frac{1}{2(1+\epsilon)}-\epsilon}\|(\phi,\psi)\|^2d\tau\leq C\delta \|(\phi,\psi)\|^2.
\end{aligned}
\end{equation}
From \cref{le3.3}, since $Z_1(t)\leq X_1(t)\leq Y_1(t),$ we have

\begin{equation}
\begin{aligned}
&C\delta\int_0^t(1+\tau)^{-\frac{1}{2}}\int_{X_1(\tau)-\sqrt{b}\tau}^{Y_1(\tau)-\sqrt{b}\tau}(\phi^2+\psi^2)e^{-\frac{\eta^2}{2\mu(1+\tau)}}d\eta d\tau\\
\leq &C\delta\int_0^t(1+\tau)^{-\frac{1}{2}}\int_{Z_1(\tau)-\sqrt{b}\tau}^{X_1(\tau)-\sqrt{b}\tau}(\phi^2+\psi^2)e^{-\frac{\eta^2}{2\mu(1+\tau)}}d\eta d\tau\\
\leq &C\delta\int_0^t(1+\tau)^{-\frac{1}{2}} e^{-\frac{(Z_1(\tau)-\sqrt{b}\tau)^2}{2\mu(1+\tau)}}\|(\phi,\psi)\|_{\infty}^2d\tau\\
\leq &C\delta \int_0^t(1+\tau)^{-\frac{1}{2}-\beta}\|(\phi,\psi)\|\|(\phi_{\eta},\psi_{\eta})\|d\tau\leq C\delta\bigg(\|(\phi,\psi)\|^2+\int_0^t\|(\phi_{\eta},\psi_{\eta})\|^2d\tau\bigg),
\end{aligned}
\end{equation}
where we use the inequality \eqref{Z1-2}.

As for the last term in \eqref{bh-s2},

\begin{equation}\label{4.64}
\begin{aligned}
&\int_0^t\int_{\mathbb{R}}\bigg[-\sqrt{b}Q_1-\mu\psi_{\eta\eta}+H\bigg](\sqrt{b}\phi-\psi)g^2 d\eta d\tau\\
\leq &C\int_0^t\int_{\mathbb{R}}|Q_1(\phi,\psi)|d\eta d\tau+C\int_0^t\|(\phi_{\eta},\psi_{\eta})\|^2d\tau\\
&+\frac{1}{160}\int_0^t\int_{-\infty}^{+\infty}w^2(\phi^2+\psi^2)d\eta d\tau+C\int_0^t\int_{\mathbb{R}}|H||(\phi,\psi)|d\eta d\tau,
\end{aligned}
\end{equation}
Inserting \eqref{J}-\eqref{4.64} into \eqref{bh-s2}, we could get

\begin{equation}\label{4.65}
\begin{aligned}
&\int_0^t\bigg<(\sqrt{b}\phi-\psi)_{\tau},(\sqrt{b}\phi-\psi)g^2\bigg>d\tau\\
\leq &C\int_0^t\|(\phi_{\eta},\psi_{\eta})\|^2d\tau+\frac{1}{80}\int_0^t\int_{-\infty}^{+\infty}w^2(\eta,\tau)(\phi^2+\psi^2)d\eta d\tau\\
&+C(\|(\phi,\psi)\|^2+1)+CG(t)+C(\delta+N(T))\int_0^t\int_{\mathbb{R}}w\bigg(\sigma(v)-\sigma(\hat{v})-\sigma'(\hat{v})\phi \bigg) d\eta d\tau.
\end{aligned}
\end{equation}

\noindent{Step 2:} 
$\eqref{pbs}_1\times\sqrt{b}+\eqref{pbs}_2$ gives that

\begin{equation}\label{4.66}
\begin{aligned}
&(\sqrt{b}\phi+\psi)_{t}-\sqrt{b}(\sqrt{b}\phi+\psi)_{\eta}-\sqrt{b}\psi_{\eta}-(\sigma(v)-\sigma(\hat{v}))_{\eta}
=-\sqrt{b}Q_1+\mu\psi_{\eta\eta}-H,
\end{aligned}
\end{equation}
and $\eqref{4.67}\times(\sqrt{b}\phi+\psi)f$ gives that

\begin{equation}\label{4.67}
\begin{aligned}
&\frac{1}{2}((\sqrt{b}\phi+\psi)^2f)_t-\frac{1}{2}(\sqrt{b}\phi+\psi)^2f_t-\frac{\sqrt{b}}{2}((\sqrt{b}\phi+\psi)^2f)_{\eta}+\frac{\sqrt{b}}{2}(\sqrt{b}\phi+\psi)^2f_{\eta}\\
&-(\sqrt{b}\phi+\psi)f(\sqrt{b}\psi+b\phi)_{\eta}-(\sqrt{b}\phi+\psi)f\bigg(\sigma(v)-\sigma(\hat{v})-b\phi\bigg)_{\eta}\\
=&(\sqrt{b}\phi+\psi)f\bigg(-\sqrt{b}Q_1+\mu\psi_{\eta\eta}-H\bigg).
\end{aligned}
\end{equation}
Note that 
\begin{equation}\label{4.68}
\begin{aligned}
\|f\|_{L^{\infty}}\leq C(1+t)^{-\frac{1}{2}},\quad \|f_t\|_{L^{\infty}}\leq C(1+t)^{-\frac{3}{2}},\quad \|f_{\eta}\|_{L^{\infty}}\leq C(1+t)^{-1}.
\end{aligned}
\end{equation}
Integrating \eqref{4.67} over $[0,t]\times\mathbb{R}$ yields 

\begin{equation}\label{4.69}
\begin{aligned}
&\sup_{t\in[0,t]}\int_{\mathbb{R}}(\sqrt{b}\phi+\psi)^2f(\eta,t)d\eta+\int_0^t\int_{\mathbb{R}}(\sqrt{b}\phi+\psi)^2w^2(\eta,t)d\eta d\tau\\
\leq &\int_{\mathbb{R}}(\sqrt{b}\phi_0+\psi_0)^2f(\eta,0)d\eta+\int_0^t\int_{\mathbb{R}}|(\phi,\psi)|^2|f_t|d\eta d\tau+\int_0^t\int_{\mathbb{R}}|(\phi,\psi)||f||Q_1,H|d\eta d\tau\\
&+\int_0^t\int_{\mathbb{R}}(\sqrt{b}\phi+\psi)f\bigg(\sigma(v)-\sigma(\hat{v})-b\phi+\mu\psi_{\eta}\bigg)_{\eta}d\eta d\tau,
\end{aligned}
\end{equation}

\begin{equation}\label{4.70}
\begin{aligned}
&\int_0^t\int_{\mathbb{R}}|(\phi,\psi)|^2|f_t|d\eta d\tau+\int_0^t\int_{\mathbb{R}}|(\phi,\psi)||f||Q_1,H|d\eta d\tau\\
\leq &C\int_0^t(1+\tau)^{-\frac{3}{2}}d\tau\int_{\mathbb{R}}|(\phi,\psi)|^2d\eta
+C\int_0^t(1+\tau)^{-2}d\tau\int_{\mathbb{R}}|(\phi,\psi)|e^{-\frac{c\eta^2}{(1+\tau)}}d\eta \\
&+C\int_0^t\|H\|_{L^1}\|(\phi,\psi)\|_{L^{\infty}}\|f\|_{L^{\infty}} d\tau,
\end{aligned}
\end{equation}

\begin{equation}\label{4.71}
\begin{aligned}
&\int_0^t\int_{\mathbb{R}}(\sqrt{b}\phi+\psi)f\bigg(\sigma(v)-\sigma(\hat{v})-b\phi+\mu\psi_{\eta}\bigg)_{\eta}d\eta d\tau\\
\leq &-\int_0^t\int_{\mathbb{R}}((\sqrt{b}\phi+\psi)f)_{\eta}\bigg(\sigma(v)-\sigma(\hat{v})-b\phi+\mu\psi_{\eta}\bigg) d\eta d\tau\\
\leq &\int_0^t\int_{\mathbb{R}}|(\phi,\psi)|w^2\bigg(\sigma(v)-\sigma(\hat{v})-b\phi\bigg) d\eta d\tau+\int_0^t\int_{\mathbb{R}}((\sqrt{b}\phi+\psi)f)_{\eta}\psi_{\eta}d\eta d\tau\\
&+\int_0^t\int_{\mathbb{R}}(\sqrt{b}\phi_{\eta}+\psi_{\eta})f\bigg(\sigma(v)-\sigma(\hat{v})-b\phi\bigg)d\eta d\tau,
\end{aligned}
\end{equation}
where the estimations for the terms in \eqref{4.71} are similar as \eqref{J}, we omit the details. Combining \eqref{4.69}-\eqref{4.71}, we have

\begin{equation}
\begin{aligned}
&\int_0^t\int_{\mathbb{R}}(\sqrt{b}\phi+\psi)^2w^2(\eta,t)d\eta d\tau\\
\leq &C\int_0^t\|(\phi_{\eta},\psi_{\eta})\|^2d\tau+\frac{1}{80}\int_0^t\int_{-\infty}^{+\infty}w^2(\eta,\tau)(\phi^2+\psi^2)d\eta d\tau\\
&+C(\|(\phi,\psi)\|^2+1)+CG(t)+C(\delta+N(T))\int_0^t\int_{\mathbb{R}}w\bigg(\sigma(v)-\sigma(\hat{v})-\sigma'(\hat{v})\phi \bigg) d\eta d\tau.
\end{aligned}
\end{equation}
Therefore, we could get \eqref{hke1-1}-\eqref{hke1-2}.

Finally, combining \eqref{hke1-1} and \eqref{hke1-2}, it yields that 

\begin{equation}
\begin{aligned}
&\int_0^t\int_{\mathbb{R}}w^2(\eta,\tau)(\phi^2+\psi^2)d\eta d\tau\\
\leq &C\int_0^t\|(\phi_{\eta},\psi_{\eta})\|^2d\tau+C(\|(\phi,\psi)\|^2+1)+CG(t)\\
&+C(\delta+N(T))\int_0^t\int_{\mathbb{R}}w\bigg(\sigma(v)-\sigma(\hat{v})-\sigma'(\hat{v})\phi \bigg) d\eta d\tau.
\end{aligned}
\end{equation}
That is \eqref{hke1}, then we complete the proof of \eqref{hke}.

\end{proof}

From \cref{4.1}-\cref{4.2}, since  
\begin{equation}\label{4.74}
\begin{aligned}
|a-v_-|\int_0^t\int_{\mathbb{R}}w\bigg(\sigma(v)-\sigma(\hat{v})-\sigma'(\hat{v})\phi \bigg) d\eta d\tau\leq CG(t),
\end{aligned}
\end{equation}
then we have 
\begin{equation}\label{4.75}
\begin{aligned}
\int_0^t\int_{\mathbb{R}}\Xi_{2x}^2(\phi^2+\psi^2)dxd\tau\leq &O(|a-v_-|^2)\int_0^t\int_{\mathbb{R}}w^2(\phi^2+\psi^2)dxd\tau\\
\leq &O(\delta)\int_0^t\|(\phi_{x},\psi_{x})\|^2d\tau+O(\delta)(\|(\phi,\psi)\|^2+1)+C\delta G(t).
\end{aligned}
\end{equation}
Inserting \eqref{4.75} into \eqref{ene1}, it yields that
\begin{equation}\label{4.76}
\begin{aligned}
\sup_{t\in[0,T]}\|(\phi,\psi)(\cdot,t)\|^2+G(t)+\int_0^t\|\psi_x\|^2dxd\tau
\leq C(\|(\phi_0,\psi_0)\|^2+1)+\frac{1}{200}\int_0^t\|\phi_x\|^2d\tau.
\end{aligned}
\end{equation}

\begin{Lem}\label{4.3}Under assumptions of \eqref{t1} in Theorem \ref{mt}, if  $N(T)\leq\epsilon_0,$ then the following energy estimate holds.

\begin{equation}\label{4.17}
\begin{aligned}
\sup_{t\in[0,T]}\|(\phi_x,\psi_x)(t)\|^2+\int_0^t\|(\phi_x,\psi_{xx})(\tau)\|^2d\tau
\leq C(\|(\phi_0,\psi_0)\|_{H^1}^2+1).
\end{aligned}
\end{equation}
\end{Lem}

\begin{proof}

$\eqref{phi-psi}_2\times\phi_x$ gives that

\begin{equation}\label{de1}
\begin{aligned}
&(\psi\phi_x)_t-(\psi\phi_t)_x+\psi_x^2-\psi_xQ_1\\
=&(\sigma(v)-\sigma(\hat{v}))_x\phi_x+(\frac{\mu}{2}\phi_x^2)_t
+\mu Q_{1x}\phi_x-H\phi_x.
\end{aligned}
\end{equation}
Integrating \eqref{de1} over $[0,t]\times\mathbb{R}$ yields 

\begin{equation}\label{dee1}
\begin{aligned}
&\sup_{t\in[0,T]}\|\phi_x(t)\|^2+\int_0^t\int_{\mathbb{R}}\sigma'(v)\phi_x^2dxd\tau\\
\leq &\|(\phi_0,\psi_0)\|_{H^1}^2+C(\|\psi(t)\|^2+\int_0^t\|\psi_x(\tau)\|^2d\tau)+C\int_0^t\int_{\mathbb{R}}\hat{v}_x^2(\sigma'(\hat{v}+\phi)-\sigma'(\hat{v}))^2 dxd\tau\\
&+\int_0^t\|(Q_1,H)\|^2d\tau.
\end{aligned}
\end{equation}
Obviously, 

\begin{equation}\label{dee1-2}
\begin{aligned}
&\int_0^t\int_{\mathbb{R}}\hat{v}_x^2(\sigma'(\hat{v}+\phi)-\sigma'(\hat{v}))^2 dxd\tau\\&\leq \int_0^t\int_{\hat{v}+\phi>a,\hat{v}>a}\hat{v}_x^2(\sigma'(\hat{v}+\phi)-\sigma'(\hat{v}))^2dxd\tau+\int_0^t\int_{\hat{v}+\phi>a,-a<\hat{v}\leq a}\hat{v}_x^2(\sigma'(\hat{v}+\phi)-\sigma'(a))^2dxd\tau\\
&~~+\int_0^t\int_{-a<\hat{v}+\phi\leq a,\hat{v}>a}\hat{v}_x^2(\sigma'(a)-\sigma'(\hat{v}))^2dxd\tau+\int_0^t\int_{-a<\hat{v}+\phi\leq a,-a<\hat{v}\leq a}\hat{v}_x^2(\sigma'(a)-\sigma'(a))^2dxd\tau\\
&\leq
CG(t) +C\delta\int_0^t\int_{-\infty}^{+\infty}w^2\phi^2dxd\tau,\\
&\int_0^t\|(Q_1,H)\|^2d\tau\leq C.
\end{aligned}
\end{equation}

$\eqref{phi-psi}_2\times(-\psi_{xx})$ gives that

\begin{equation}\label{de2}
\begin{aligned}
-(\psi_t\psi_x)_x+(\frac{1}{2}\psi_x^2)_t+\mu\psi_{xx}^2+\bigg\{\sigma'(v)(\hat{v}_x+\phi_x)-\sigma'(\hat{v})\hat{v}_x\bigg\}\psi_{xx}=H\psi_{xx}.
\end{aligned}
\end{equation}
Integrating \eqref{de2} over $[0,t]\times\mathbb{R}$ yields that

\begin{equation}\label{dee2}
\begin{aligned}
&\sup_{t\in[0,T]}\|\psi_x(t)\|^2+\int_0^t\|\psi_{xx}(\tau)\|^2d\tau\\
\leq &\|(\phi_0,\psi_0)\|_{H^1}^2+C\int_0^t\|\phi_x(\tau)\|^2d\tau+\int_0^t\int_{\mathbb{R}}\hat{v}_x^2(\sigma'(\hat{v}+\phi)-\sigma'(\hat{v}))^2 dxd\tau\\
&+C\int_0^t\|H(\tau)\|^2d\tau.
\end{aligned}
\end{equation}
Combining \eqref{dee1},\eqref{dee1-2} ,\eqref{dee2} and previous results, let $\delta$ is suitably small, we could get \eqref{4.17}.

\end{proof}

Now we turn to show $(\phi,\psi)(x,t)\rightarrow 0$ as $t\rightarrow \infty.$ Combining lemma \cref{4.1}-\cref{4.3}, one can extend the local solution globally and show that the following a priori estimates holds for $t\in[0,+\infty).$

\begin{equation}
\begin{aligned}
\sup_{t\in[0,T]}\|(\phi,\psi)(t)\|_{H^1}^2+\int_0^t\|(\phi_x,\psi_x,\psi_{xx})(\tau)\|^2d\tau+G(t)
\leq C(\|(\phi_0,\psi_0)\|_{H^1}^2+1).
\end{aligned}
\end{equation}

It tell us that

\begin{equation}
\begin{aligned}
\int_0^{+\infty}\bigg\{\|(\phi_{x},\psi_{x})(\cdot,t)\|^2+\bigg|\frac{\partial}{\partial t}\|(\phi_{x},\psi_{x})(\cdot,t)\|^2\bigg|\bigg\}dt<+\infty,
\end{aligned}
\end{equation}
this gives that 

\begin{equation}
\begin{aligned}
&\lim_{t\rightarrow+\infty}\|(\phi_{x},\psi_{x})(\cdot,t)\|^2=0,\\
&\lim_{t\rightarrow+\infty}\|(\phi,\psi)(\cdot,t)\|_{L^{\infty}(\mathbb{R})}^2\leq  \lim_{t\rightarrow+\infty}\bigg\{2\|\phi(\cdot,t)\|\|\phi_{x}(\cdot,t)\|+2\|\psi(\cdot,t)\|\|\psi_{x}(\cdot,t)\|\bigg\}=0,
\end{aligned}
\end{equation}
we complete the proof of theorem \ref{mt}.

\section*{Acknowledgements}
Guo is partially supported by the National Natural Science Foundation of China under contract No.11931013 and the Natural Science Foundation of Guangxi Province under contract No. 2022GXNSFDA035078. Hou is supported by the Scientific Research Program Funded by Education Department of Shaanxi Provincial Government (Program No.22JK0582). Lingda Xu is supported by the Research Centre for Nonlinear Analysis at The Hong Kong Polytechnic University.

\section*{Data Availability Statement}

Data sharing is not applicable to this article as no data sets were generated or analysed during the current study.

\end{document}